\documentclass[final,11pt]{elsarticle}

\usepackage{verbatim,a4wide}

\usepackage{amsmath}
\usepackage{amssymb}
\usepackage{amsthm}
\usepackage{bm}

\usepackage{algorithm}
\usepackage{algpseudocode}

\usepackage{xcolor}

\usepackage{tikz}
\usetikzlibrary{patterns}

\usepackage{hyperref}

\usepackage[cp1250]{inputenc}
\usepackage[OT4]{fontenc}
\usepackage[english]{babel}

\numberwithin{equation}{section}
\numberwithin{algorithm}{section}
\numberwithin{figure}{section}
\numberwithin{table}{section}

\newtheorem{theorem}{Theorem}[section]

\newtheorem{remark}[theorem]{Remark}

\newtheorem{example}[theorem]{Example}

\newtheorem{problem}[theorem]{Problem}

\newcommand{\p}[1]{\mbox{\textsf{#1}}}

\journal{Computer Aided Geometric Design} 

\begin{document}

%%%%%%%%%%%%%%%%%%%%%%%%%%%%%%%%%%%%%%%%%%%%%%%%%%%%%%%%%%%%%%%%%%%%%%%%%%%%%%
%%%%%%%%%%%%%%%%%%%%%%%%%%%%%%%%%%%%%%%%%%%%%%%%%%%%%%%%%%%%%%%%%%%%%%%%%%%%%%
%%%%%%%%%%%%%%%%%%%%%%%%%%%%%%%%%%%%%%%%%%%%%%%%%%%%%%%%%%%%%%%%%%%%%%%%%%%%%%
%%%%%%%%%%%%%%%%%%%%%%%%%%%%%%%%%%%%%%%%%%%%%%%%%%%%%%%%%%%%%%%%%%%%%%%%%%%%%%
\begin{frontmatter}

%%%%%%%%%%
\title{Fast evaluation of derivatives of B\'{e}zier curves}
%%%%%%%%%%

\author[A1,A2]{Filip Chudy\corref{cor}}
\ead{Filip.Chudy@cs.uni.wroc.pl}

\author[A1]{\texorpdfstring{Pawe{\l} Wo\'{z}ny}{Pawel Wozny}}
\ead{Pawel.Wozny@cs.uni.wroc.pl}

\cortext[cor]{Corresponding author.}

\address[A1]{Institute of Computer Science, University of Wroc{\l}aw,
             ul.~Joliot-Curie 15, 50-383 Wroc{\l}aw, Poland}         

\address[A2]{Department of Mathematics, University of Bologna,
             Piazza di Porta San Donato 5, 40126 Bologna, Italy}

\begin{abstract}
New geometric methods for fast evaluation of derivatives of polynomial and 
rational B\'{e}zier curves are proposed. They apply an algorithm for evaluating 
polynomial or rational B\'{e}zier curves, which was recently given by the 
authors. Numerical tests show that the new approach is more 
efficient than the methods which use the famous de Casteljau algorithm. The 
algorithms work well even for high-order derivatives of rational B\'{e}zier 
curves of high degrees. 
\end{abstract}

%%%%%%%%%%%%%%%%%%%%%%%%%%%%%%%%%%%%%%%%%%%%%%%%%%%%%%%%%%%%%%%%%%%%%%%%%%%%%%
% ToDo
%
%\begin{highlights}
%
%\item ToDo
%
%
%\end{highlights}
%%%%%%%%%%%%%%%%%%%%%%%%%%%%%%%%%%%%%%%%%%%%%%%%%%%%%%%%%%%%%%%%%%%%%%%%%%%%%%


\end{frontmatter}
%%%%%%%%%%%%%%%%%%%%%%%%%%%%%%%%%%%%%%%%%%%%%%%%%%%%%%%%%%%%%%%%%%%%%%%%%%%%%%
%%%%%%%%%%%%%%%%%%%%%%%%%%%%%%%%%%%%%%%%%%%%%%%%%%%%%%%%%%%%%%%%%%%%%%%%%%%%%%
%%%%%%%%%%%%%%%%%%%%%%%%%%%%%%%%%%%%%%%%%%%%%%%%%%%%%%%%%%%%%%%%%%%%%%%%%%%%%%
%%%%%%%%%%%%%%%%%%%%%%%%%%%%%%%%%%%%%%%%%%%%%%%%%%%%%%%%%%%%%%%%%%%%%%%%%%%%%%

%%%%%%%%%%%%%%%%%%%%%%%%%%%%%%%%%%%%%%%%%%%%%%%%%%%%%%%%%%%%%%%%%%%%%%%%%%%%%%
%%%%%%%%%%%%%%%%%%%%%%%%%%%%%%%%%%%%%%%%%%%%%%%%%%%%%%%%%%%%%%%%%%%%%%%%%%%%%%
\section{Introduction}                                  \label{S:Introduction}
%%%%%%%%%%%%%%%%%%%%%%%%%%%%%%%%%%%%%%%%%%%%%%%%%%%%%%%%%%%%%%%%%%%%%%%%%%%%%%
%%%%%%%%%%%%%%%%%%%%%%%%%%%%%%%%%%%%%%%%%%%%%%%%%%%%%%%%%%%%%%%%%%%%%%%%%%%%%%

The $k$th Bernstein polynomial of degree $n$ is given by the formula
$$
B^n_k(t):=\binom{n}{k}t^k(1-t)^{n-k} 
                    \qquad (k=0,1,\ldots,n;\, n\in\mathbb{N}).
$$
In the sequel, a convention is applied that $B^n_k\equiv0$ if $k<0$ or 
$k>n$.

The polynomials $B^n_0,B^n_1,\ldots,B^n_n$ form a basis of the space of all 
polynomials of degree at most~$n$. They satisfy the degree reduction and degree 
elevation relations:
\begin{eqnarray}\label{E:BernDegDown}
B^n_k(t) & = & t B^{n-1}_{k-1}(t)+ (1-t) B^{n-1}_k(t),\\
\label{E:BernDegUp} 
B^n_k(t) & = & \dfrac{n-k+1}{n+1} B^{n+1}_{k}(t)
                              + \dfrac{k+1}{n+1} B^{n+1}_{k+1}(t),
\end{eqnarray}
where $k=0,1,\ldots,n$. The following identity is also true:
\begin{equation}\label{E:BernsteinDerivSimple}
\left(B^n_k(t)\right)' = n \left(B^{n-1}_{k-1}(t) - B^{n-1}_k(t) \right)
                                                 \qquad (0\leq k\leq n).
\end{equation}
See, e.g.,~\cite{Farin2002}.

Due to their properties, such as the partition of unity for $t \in \mathbb{R}$
and non-negativity for $t \in [0, 1]$, Bernstein polynomials have found many
applications in computer-aided design, approximation theory and numerical
analysis. For more properties, history and applications of Bernstein 
polynomials, see, e.g., \cite{Farouki2012}.

A rational B\'{e}zier curve $\p{R}_n : [0, 1]\rightarrow\mathbb{E}^d$ of degree
$n$ with control points $\p{W}_0, \p{W}_1, \ldots, \p{W}_n \in \mathbb{E}^d$
and their corresponding weights $\omega_0, \omega_1, \ldots, \omega_n > 0$
is defined as
\begin{equation}\label{E:RationalBezierCurve}
\p{R}_n(t):=\frac{\displaystyle \sum_{k=0}^n\omega_k B^n_k(t)\p{W}_k} %\Big/ 
                 {\displaystyle \sum_{k=0}^n\omega_k B^n_k(t)}
                                                      \qquad (0\leq t\leq 1).
\end{equation}
 
If $\omega_0 = \omega_1 = \ldots = \omega_n$, a rational B\'{e}zier curve 
reduces to a polynomial one. For clarity, a~polynomial B\'{e}zier curve of 
degree $n$ will be denoted as $\p{P}_n$.

Rational B\'{e}zier curves are among the most fundamental tools in computer 
graphics. The classic method of evaluating these curves is by using the famous 
de Casteljau algorithm (see, e.g., \cite{BM99,DC59,DC63,DC99,Farin2002}) which 
is the consequence of 
the relation~\eqref{E:BernDegDown}:
\begin{eqnarray}\nonumber
&&\p{W}_k^{(0)}:=\p{W}_k,\qquad \omega_k^{(0)}:=\omega_k
                                            \qquad (0\leq k\leq n),\\[1ex]
&&\label{E:RatdeCastel}
\left.
\begin{array}{l}
\omega^{(i)}_k:=(1-t)\omega^{(i-1)}_k+t\omega^{(i-1)}_{k+1},\\[1ex]
\p{W}^{(i)}_k:=(1-t)\dfrac{\omega^{(i-1)}_k}{\omega^{(i)}_k}\p{W}^{(i-1)}_k
+t\dfrac{\omega^{(i-1)}_{k+1}}{\omega^{(i)}_k}\p{W}^{(i-1)}_{k+1}\\[1ex]
\end{array}
\right\}
\qquad (1\leq i\leq n;\; 0\leq k\leq n-i).\qquad
\end{eqnarray}
Then, $\p{R}_n(t)=\p{W}^{(n)}_0$. For polynomial B\'{e}zier curves, the 
algorithm is simpler, as the quantities $\omega^{(i)}_k$ can be omitted.

The method has a geometric interpretation, good numerical properties, and 
computes only convex combinations of points in ${\mathbb E^d}$. This algorithm 
can also be used for subdividing a B\'{e}zier curve into two parts. In 
a polynomial case, it (almost) computes derivatives of a~B\'{e}zier curve. 
However, the main drawback of the de Casteljau algorithm is its complexity, as 
it requires $O(n^2 d)$ operations to evaluate a point on a polynomial or 
rational B\'{e}zier curve of degree $n$ in $\mathbb{E}^d$.

Recently, the authors have proposed a new linear-time geometric methods for 
evaluation polynomial and rational B\'{e}zier curves which retain the neat
geometric and numerical qualities of the de Casteljau algorithm, while reducing
the complexity to $O(nd)$. Numerical experiments have shown that the new
algorithms are much faster than these using the de Casteljau algorithm. For
details, efficient implementations, as well as results of numerical tests, see
\cite{WCh2020} and \cite[Chapter 2]{FChPhD}.

Now, we briefly present the main idea of this method focusing on its 
\textit{geometricity} which has not been fully done in~\cite{WCh2020}, where the 
method was described for the first time. We also recall the most important 
properties of this approach for fast evaluation of B\'{e}zier curves.   

One can reformulate Eq.~\eqref{E:RationalBezierCurve} as
\begin{eqnarray}
\nonumber
\p{R}_n(t)&=&\dfrac{\displaystyle \sum_{k=0}^i \omega_k B^n_k(t)}
                   {\displaystyle \sum_{k=0}^n \omega_k B^n_k(t)} \cdot
             \dfrac{\displaystyle \sum_{k=0}^i \omega_k B^n_k(t) \p{W}_k}
                   {\displaystyle \sum_{k=0}^i \omega_k B^n_k(t)}
           + \dfrac{\displaystyle \sum_{k=i+1}^n \omega_k B^n_k(t) \p{W}_k}
                   {\displaystyle \sum_{k=0}^n \omega_k B^n_k(t)}\\
\label{E:RationalBezierPreQ}                   
&=&
\dfrac{\displaystyle \sum_{k=0}^i \omega_k B^n_k(t)}
      {\displaystyle \sum_{k=0}^n \omega_k B^n_k(t)} \cdot \p{Q}^n_i(t) +
\dfrac{\displaystyle \sum_{k=i+1}^n \omega_k B^n_k(t) \p{W}_k}
      {\displaystyle \sum_{k=0}^n \omega_k B^n_k(t)},
\end{eqnarray}
where 
\begin{equation}\label{E:Def_Q}
 \p{Q}_i\equiv \p{Q}^n_i(t):= \sum_{k=0}^i \omega_k B^n_k(t) \p{W}_k \Big/
                              \sum_{k=0}^i \omega_k B^n_k(t)
\end{equation}
and $i=0,1,\ldots,n$. Note that $\p{Q}_i$ $(0\leq i\leq n)$ is a well-defined 
point in $\mathbb{E}^d$, as it is a barycentric (and convex) combination of 
points in $\mathbb{E}^d$. In particular, $\p{Q}_0 = \p{W}_0$ and, for $i=n$, 
$\p{R}_n(t) = \p{Q}_n$.

Comparing the right-hand sides of Eq.~\eqref{E:RationalBezierPreQ} for $i$ and
$i-1$ gives, after simple algebra, the relation
$$
\p{Q}_i =
\dfrac{\displaystyle \sum_{k=0}^{i-1} \omega_k B^n_k(t)}
      {\displaystyle \sum_{k=0}^i \omega_k B^n_k(t)} \cdot
\p{Q}_{i-1}
+
\dfrac{\omega_i B^n_i(t)}
      {\displaystyle \sum_{k=0}^i\omega_k B^n_k(t)} \cdot
 \p{W}_i.
$$
Let
\begin{equation}\label{E:BezierDefHSum}
h_i\equiv h^n_i(t):=\omega_i B^n_i(t)\Big/\sum_{k=0}^i\omega_k B^n_k(t)
                                                     \qquad (i=1,2,\ldots,n),
\end{equation}
with $h^n_0(t):=1$.

From~\eqref{E:BezierDefHSum}, it follows that
$$
  h^n_i(t) = 
      \dfrac{t (n-i+1) \omega_{i} h^n_{i-1}(t)}
            {(1-t) i \omega_{i-1} + t (n-i+1) \omega_{i} h^n_{i-1}(t)} 
                                                     \qquad (i=1,2,\ldots,n).
$$

To sum up, we get the following recurrence scheme:
\begin{equation}\label{E:Def_h_k_Q_k} 
\left\{
 \begin{array}{l}
 h_0:=1,\quad \p{Q}_0:=\p{W}_0,\\[1ex]
 h_i:=\dfrac{\displaystyle \omega_ih_{i-1}t(n-i+1)}
            {\displaystyle \omega_{i-1}i(1-t)+\omega_ih_{i-1}t(n-i+1)},\\[2.5ex]
 \p{Q}_i:=(1-h_i)\p{Q}_{i-1}+h_i\p{W}_i,
\end{array}
\right.
\end{equation}
where $1 \leq i \leq n$. For $t \in [0, 1]$, the main properties of the 
quantities evaluated by the scheme~\eqref{E:Def_h_k_Q_k} are as follows:
\begin{enumerate}
  \itemsep1ex
  \item $h_i\in[0,1]$ (thus $\p{Q}_i$ is a convex combination of $\p{Q}_{i-1}$
  	    and $\p{W}_i$),
  \item $\p{Q}_i\in\mathbb E^d$,
  \item $\p{Q}_i\in C_i \equiv \mbox{conv}\{\p{W}_0,\p{W}_1,\ldots,\p{W}_i\}$
          (thus $\mbox{conv}\{\p{Q}_0,\p{Q}_1,\ldots,\p{Q}_i\}\subseteq C_i$),
  \item $\p{R}_n([0,u])\subseteq \mbox{conv}\{\p{Q}_0, \p{Q}_1, \ldots,
         \p{Q}_n\}\quad (u\leq t)$,
\end{enumerate}
where $\mbox{conv} A$ is the convex hull of a set $A$. Moreover,
$\p{R}_n(t)=\p{Q}_n$. 

\begin{figure}[ht!]
  \centering
  \vspace*{-3.2ex}
  \includegraphics[width=\textwidth]{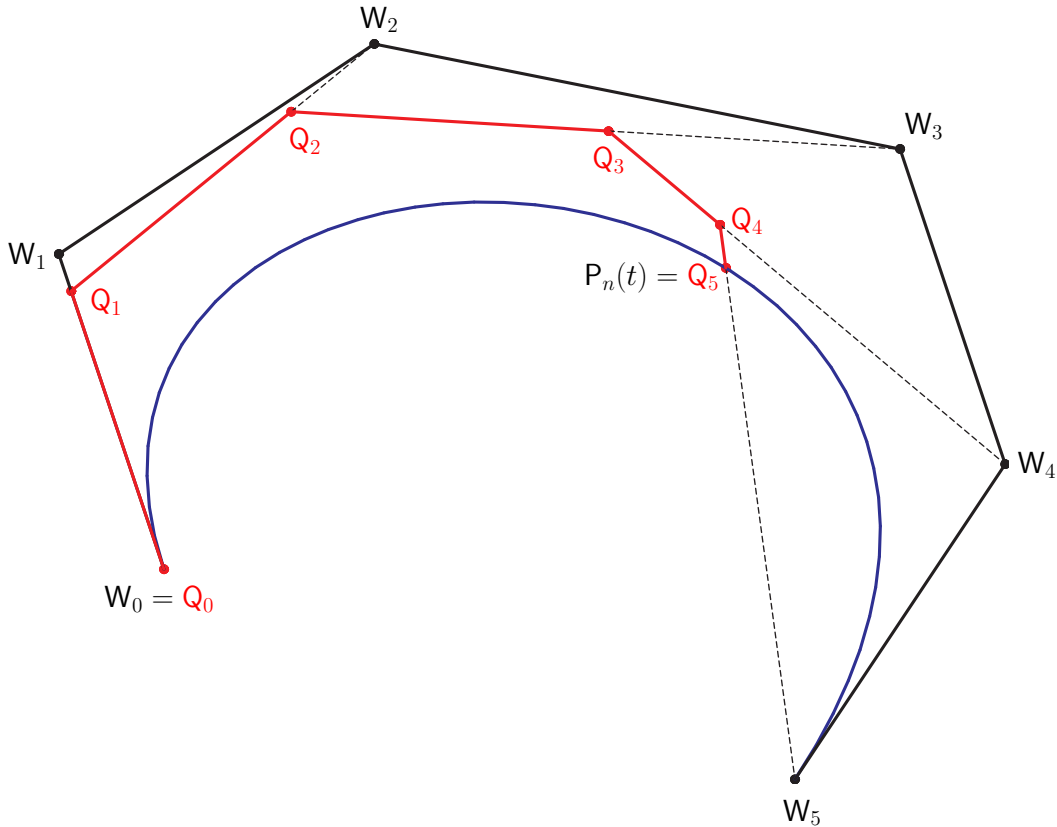}
  \caption{A computation of a point on a planar polynomial B\'{e}zier curve of  
           degree $n=5$ using the new method. The method computes much fewer
           intermediate points compared to the de Casteljau algorithm. Image
           taken from~\cite{WCh2020}.}\label{F:Figure1}
\end{figure}

From the given properties, it follows that the new algorithm has geometric 
interpretation (see Figure~\ref{F:Figure1}), computes only convex combination 
of points, has the \textit{convex hull property} and its computational 
complexity is linear with respect to the number of control points, i.e., it is 
asymptotically optimal. It was also shown that the proposed method can find 
application in B\'{e}zier curves subdivision. 

Certainly, the method simplifies for polynomial B\'{e}zier curves. 
Figure~\ref{F:Figure1} illustrates the new method in the case of a~planar 
polynomial B\'{e}zier curve of degree $n=5$. 

\begin{remark}\label{R:ManyPolyCurves}
It is worth mentioning that the new algorithm is especially useful if one has to
evaluate many polynomial B\'{e}zier curves of the same degree and for the same
parameter $t\in[0,1]$ --- as, for example, in the problem of rendering of
rectangular polynomial B\'{e}zier patches --- because the quantities $h^n_i(t)$ 
do not depend on the control points and thus can be computed just once for all 
curves. This fact will also find its application when evaluating the derivatives 
of polynomial and rational B\'{e}zier curves.
\end{remark}

The new algorithm can be used, as seen in \cite{RH21}, to pre-compute some
quantities to further accelerate the computations for B\'{e}zier curves of very
high degrees or for computing many points on a single curve --- however, at the
cost of losing the geometric interpretation.

Additionally, in~\cite{RV21}, it is shown how to adapt the new algorithm to
evaluate algebraic-hyperbolic PH curves (EPH curves) for a fixed parameter $t$.
The approach presented there is to convert an EPH curve into a polynomial
B\'{e}zier curve such that they have the same value at~$t$, and then to evaluate
the newly found B\'{e}zier curve using the method presented above.

Recently, this method has also been used in~\cite{ChW2023} (see also 
\cite[Chapter 3]{FChPhD}), where a fast method of evaluating B-spline curves 
was proposed. 

It is possible to generalize the new approach and use it for a broader family of 
\textit{rational parametric objects} (in particular, for rectangular and 
triangular rational B\'{e}zier surfaces; cf.~\cite[\S3]{WCh2020} and 
\cite[Chapter 2]{FChPhD}) $\p{S}_N:C\rightarrow \mathbb E^d$ $(d\in\mathbb N)$ 
of the form
\begin{equation}\label{E:ParObject}  
\p{S}_N(\bm{t}):=\frac{\displaystyle \sum_{k=0}^{N}\omega_k\p{W}_kb_k(\bm{t})}
                      {\displaystyle \sum_{k=0}^{N}\omega_k b_k(\bm{t})},
\end{equation}
with the weights $\omega_k>0$, and control points $\p{W}_k\in\mathbb E^d$ 
$(0\leq k\leq N)$, where $b_k:D\rightarrow\mathbb R$ 
$(k=0,1,\ldots,N;\, N\in\mathbb N)$ are the real-valued multivariable 
\textit{basis functions} such that
\begin{equation}\label{E:BasisF}
b_k(\bm{t})\geq0,\quad \sum_{k=0}^{N}b_k(\bm{t})\equiv 1
\end{equation} 
for $\bm{t}\in C\subseteq D$. See \cite[\S1 and Algorithm~1.1]{WCh2020}. 

The main goal of this paper is to show that the new geometric method for fast
evaluation of rational B\'{e}zier curves can also be applied to efficiently 
comput the derivatives of such curves in a \textit{geometric way}. We consider
the following problem.

\begin{problem}\label{P:ComputeDeriv}
For $n,r\in\mathbb{N}$, control points $\p{W}_0, \p{W}_1, \ldots, \p{W}_n \in
\mathbb{E}^d$, their associated weights $\omega_0, \omega_1, \ldots, \omega_n 
\in \mathbb{R}_+$, $t \in [0, 1]$, propose a {\rm geometric} method of 
evaluating
$$
\p{R}_n(t),\p{R}_n^{(1)}(t),\ldots,\p{R}_n^{(r)}(t)
$$
(cf.~\eqref{E:RationalBezierCurve}).
\end{problem}

In Section~\ref{S:PolyDiff}, we focus on the derivatives of polynomial 
B\'{e}zier curves, and in Section~\ref{S:RatDiff} we discuss 
Problem~\ref{P:ComputeDeriv} in its general setting.

More precisely, we present five different approaches which all use the
linear-time geometric methods for evaluating polynomial and rational B\'{e}zier
curves proposed by the authors in~\cite{WCh2020}.

The first two algorithms allow to compute derivatives of polynomial B\'{e}zier 
curves. The method from~\S\ref{SS:PolyDiff-I} is very efficient if $d\geq2$
(especially for small $n,r$ or if $r=n$), while its variant given 
in~\S\ref{SS:PolyDiff-II} is much better if, for example, $d=1$ and 
$n\geq 20$, i.e., for a linear combination of Bernstein polynomials which is 
very important in many applications. Notice that both methods outperform the
de~Casteljau algorithm.

The three methods given in Section~\ref{S:RatDiff} compute derivatives of 
rational B\'{e}zier curves. The algorithm proposed in~\S\ref{SS:RatDiff-Floater} 
accelerates Floater's formulas found in~\cite{Floater91} and only applies to 
the first and the second derivative. The fourth and the fifth method allow to 
evaluate the derivatives of any order. The approach presented 
in~\S\ref{SS:RatDiff-II} is simpler, but it is slower in general. On the other 
hand, the efficient implementation of the algorithm from~\S\ref{SS:RatDiff-III} 
is more complicated, but it results in a faster method. 

In this paper, the advantages and the disadvantages of the new methods will be examined,
including their computational complexity and numerical behavior.
Due to the fact that none of the methods universally outperforms others
(or cannot be used for derivatives of any order),
the reader needs to choose the approach which is best suited for their
particular use case.

%%%%%%%%%%%%%%%%%%%%%%%%%%%%%%%%%%%%%%%%%%%%%%%%%%%%%%%%%%%%%%%%%%%%%%%%%%%%%%
%%%%%%%%%%%%%%%%%%%%%%%%%%%%%%%%%%%%%%%%%%%%%%%%%%%%%%%%%%%%%%%%%%%%%%%%%%%%%%
\section{Derivatives of polynomial B\'{e}zier curves}       \label{S:PolyDiff}
%%%%%%%%%%%%%%%%%%%%%%%%%%%%%%%%%%%%%%%%%%%%%%%%%%%%%%%%%%%%%%%%%%%%%%%%%%%%%%
%%%%%%%%%%%%%%%%%%%%%%%%%%%%%%%%%%%%%%%%%%%%%%%%%%%%%%%%%%%%%%%%%%%%%%%%%%%%%%

Assume that $t\in[0,1]$. Let $\p{P}_n$ be a polynomial B\'{e}zier curve of 
degree $n$ with control points $\p{W}_0, \p{W}_1,\ldots,\p{W}_n\in\mathbb{E}^d$:
\begin{equation}\label{E:PolyBezierCurve}
\p{P}_n(t):=\sum_{k=0}^{n}B^n_k(t)\p{W}_k. 
\end{equation}
Differentiating both sides of Eq.~\eqref{E:PolyBezierCurve} and 
applying~\eqref{E:BernsteinDerivSimple} gives
\begin{equation}\label{E:PolyBezierFirstDeriv} 
\p{P}_n'(t) = \sum_{k=0}^{n-1} B^{n-1}_{k}(t) \p{v}_k^{(1)},
\end{equation}
where the vectors $\p{v}_k^{(1)}\in{\mathbb R}^d$ ($k=0,1,\ldots,n-1$) are 
given by
\begin{equation}\label{E:ControlVectors-0}
\p{v}_k^{(1)} := n \left( \p{W}_{k+1} - \p{W}_k \right).
\end{equation}
A \textit{vector B\'{e}zier curve} of degree $n-1$ with \textit{control vectors} 
$\p{v}_k^{(1)}\in{\mathbb R}^d$ has been obtained. Repeating this process gives 
an arbitrary derivative of a polynomial B\'{e}zier curve:
$$
\p{P}_n^{(j)}(t) = \sum_{k=0}^{n-j} B^{n-j}_k(t) \p{v}_k^{(j)}
                                                        \qquad(2\leq j\leq n),
$$
where
\begin{equation}\label{E:ControlVectors}
\p{v}_k^{(j)} := (n-j+1)\left( \p{v}_{k+1}^{(j-1)} - \p{v}_k^{(j-1)} \right)
                                                    \qquad (k=0,1,\ldots,n-j).
\end{equation}

Certainly, this is a classic result. Moreover, it is also well-known that
\begin{equation}\label{E:r-DiffPolyBezier}
\p{P}_n^{(j)}(t)=\dfrac{n!}{(n-j)!}\Delta^j\p{W}^{(n-j)}_0(t)
                                                       \qquad (j\leq n),
\end{equation}
where $\p{W}^{(n-j)}_k(t)\equiv \p{W}^{(n-j)}_k$ $(0\leq k\leq j)$ are the 
points from the $(n-j)$th column of the table computed using the polynomial 
de Casteljau algorithm (cf.~\eqref{E:RatdeCastel}). Here, $\Delta$ is the 
forward difference operator, i.e.,
\begin{equation}\label{E:Delta}
\Delta c_k:=c_{k+1}-c_k.
\end{equation}
In particular, we have
\begin{eqnarray*}
&&P_n^{(j)}(0)=\dfrac{n!}{(n-j)!}\sum_{k=0}^j\binom{j}{k}(-1)^{j-k}\p{W}_k,\\
&&P_n^{(j)}(1)=\dfrac{n!}{(n-j)!}\sum_{k=n-j}^n\binom{j}{n-k}(-1)^{n-k}\p{W}_k.
\end{eqnarray*}
See, e.g., \cite[Eq.~(5.25) and (5.26)]{Farin2002}.

From~\eqref{E:r-DiffPolyBezier} and the polynomial de Casteljau algorithm, it 
follows that, for a given $t\in[0,1]$ and $0\leq r\leq n$, one can compute the 
point $\p{P}_n(t)\in{\mathbb E}^d$ and all the vectors 
$$
\p{P}_n'(t),\p{P}_n''(t),\ldots,\p{P}_n^{(r)}(t)\in{\mathbb R}^d
$$
(cf.~Problem~\ref{P:ComputeDeriv}) with total $O(d(n^2+r^2))$ 
complexity.

%%%%%%%%%%%%%%%%%%%%%%%%%%%%%%%%%%%%%%%%%%%%%%%%%%%%%%%%%%%%%%%%%%%%%%%%%%%%%%
\subsection{New method}                                  \label{SS:PolyDiff-I}
%%%%%%%%%%%%%%%%%%%%%%%%%%%%%%%%%%%%%%%%%%%%%%%%%%%%%%%%%%%%%%%%%%%%%%%%%%%%%%

Now, we show that a polynomial version of Problem~\ref{P:ComputeDeriv} can be 
solved much faster. Let us fix $0\leq r\leq n$. First, we compute all the
control vectors $\p{v}_k^{(j)}\in{\mathbb R}^d$ for $j=1,2,\ldots,r$ and 
$k=0,1,\ldots,n-j$. It can be done in $O(rnd)$ time 
(cf.~\eqref{E:ControlVectors-0}, \eqref{E:ControlVectors}). Then we apply 
the polynomial version of \cite[Algorithm 2.2]{WCh2020}---which gives a very 
efficient and numericaly stable implementation of the 
scheme~\eqref{E:Def_h_k_Q_k}---$(r+1)$ times, once for each derivative. 

Taking into account that the cost of the new method of evaluation of B\'{e}zier 
curves is linear with respect to the number of control points 
(or control vectors) (cf.~Section~\ref{S:Introduction}), the total complexity 
of this approach is $O(rnd)$. Algorithm~\ref{A:PolyBezierDeriv} implements this 
method, where \textsc{NewPolyBezierEval} denotes the polynomial version of 
\cite[Algorithm 2.2]{WCh2020}. 

\begin{algorithm}[ht!]
\caption{Algorithm for finding the value and first $r$ derivatives of 
a polynomial B\'{e}zier curve of degree $n$ at $t$.}\label{A:PolyBezierDeriv}
\begin{algorithmic}[1]
\Procedure {NewPolyBezierDerivEval}{$n, t, r, \p{W}$}

\For {$k \gets 0,n$}
    \State $\p{v}^{(0)}_k \gets \p{W}_k$
\EndFor

\For {$j \gets 1,r$}
    \For {$k \gets 0,n-j$}
        \State $\p{v}^{(j)}_k\gets 
                   (n-j+1)\cdot(\p{v}^{(j-1)}_{k+1}-\p{v}^{(j-1)}_k)$
    \EndFor
\EndFor

\For {$j \gets 0,r$}
    \State $\p{P}_{n}^{(j)}(t)\gets 
          \Call{NewPolyBezierEval}{n-j,t,\p{v}^{(j)}_0,\p{v}^{(j)}_1,
                                                    \ldots,\p{v}^{(j)}_{n-j}}$
\EndFor

\State \Return $\p{P}_{n}^{(0)}(t),\p{P}_{n}^{(1)}(t),\ldots,\p{P}_{n}^{(r)}(t)$

\EndProcedure
\end{algorithmic}
\end{algorithm}

\begin{remark}\label{R:Remark2}
Similarly to Remark~\ref{R:ManyPolyCurves}, observe that if one computes the 
$j$th derivative of many polynomial B\'{e}zier curves of the 
same degree $n$ at the same point $t$ using the presented approach then the 
quantities $h^{n-j}_i(t)$ (cf.~\eqref{E:BezierDefHSum} and 
\eqref{E:Def_h_k_Q_k}) can be evaluated just once. Then, the computation of 
$h^{n-j}_i(t)$ in every call of \textsc{NewPolyBezierEval} in
Algorithm~\ref{A:PolyBezierDeriv} can be skipped, which speeds up the whole 
process.  
\end{remark}

%%%%%%%%%%%%%%%%%%%%%%%%%%%%%%%%%%%%%%%%%%%%%%%%%%%%%%%%%%%%%%%%%%%%%%%%%%%%%%
\subsection{Modification of the new method}             \label{SS:PolyDiff-II}
%%%%%%%%%%%%%%%%%%%%%%%%%%%%%%%%%%%%%%%%%%%%%%%%%%%%%%%%%%%%%%%%%%%%%%%%%%%%%%

Certainly, to solve the polynomial version of Problem~\ref{P:ComputeDeriv}
using Algorithm~\ref{A:PolyBezierDeriv}, it is necessary to compute the 
quantities $h^{n-j}_i(t)$ for all $j=0,1,\ldots,r$ and all $i=0,1,\ldots,n-j$.

Observe that it is possible to find the point $\p{P}_{n}(t)\in{\mathbb E}^d$ 
and all the vectors 
\begin{equation}\label{E:AllPolyDeriv}
\p{P}_{n}^{(1)}(t),\p{P}_{n}^{(2)}(t),\ldots,\p{P}_{n}^{(r)}(t)\in{\mathbb R}^d
\end{equation}
using only the values $h^{n}_i(t)$ for $i=0,1,\ldots,n$. 

Indeed, in order to do this, one has to write all 
derivatives~\eqref{E:AllPolyDeriv} as linear combinations of Bernstein 
polynomials of degree $n$ and then evaluate them, taking  
Remark~\ref{R:ManyPolyCurves} into account. More precisely, from 
relations~\eqref{E:BernsteinDerivSimple} and \eqref{E:BernDegUp}, it follows 
that
\begin{equation}\label{E:BernsteinDeriv}
\left(B^n_k(t)\right)'=(n-k+1)B^n_{k-1}(t)+(2k-n)B^n_k(t)-(k+1)B^n_{k+1}(t)
                                                       \qquad (0\leq k\leq n).
\end{equation}
Using this identity, we obtain
$$   
\p{P}_n'(t)=\sum_{k=0}^{n}B^n_{k}(t)\p{u}_k^{(1)},
$$
where the vectors $\p{u}_k^{(1)}$ ($k=0,1,\ldots,n$) are given by
$$
\p{u}_k^{(1)}:=(n-k)\p{W}_{k+1}+(2k-n)\p{W}_k-k\p{W}_{k-1}
$$
(cf.~\eqref{E:PolyBezierFirstDeriv}, \eqref{E:ControlVectors-0}). It is not
difficult to see that for $j=1,2,\ldots,n$, we have
$$
\p{P}_n^{(j)}(t)=\sum_{k=0}^{n}B^n_{k}(t)\p{u}_k^{(j)},
$$
where vectors $\p{u}_k^{(j)}$ are defined recursively in the following way:
$$
\p{u}_k^{(j)}:=(n-k)\p{u}_{k+1}^{(j-1)}+(2k-n)\p{u}_{k}^{(j-1)}
                                            -k\p{u}_{k-1}^{(r-1)}=
               (n-k)\Delta \p{u}_k^{(j-1)}+k\Delta \p{u}_{k-1}^{(j-1)}
                                                       \quad (0\leq k\leq n)
$$                  
(cf.~\eqref{E:Delta}). Algorithm~\ref{A:ModPolyBezierDeriv} implements this 
method. Certainly, the presented approach's complexity is also $O(rnd)$ 
(cf.~Algorithm~\ref{A:PolyBezierDeriv}). 

Notice that in the case of this approach the improvement similar to the one 
mentioned in Remark~\ref{R:Remark2} is also possible.

\begin{algorithm}[ht!]
\caption{Modified algorithm for finding the value and first $r$ derivatives of 
a polynomial B\'{e}zier curve of degree $n$ at $t$.}\label{A:ModPolyBezierDeriv}
\begin{algorithmic}[1]
\Procedure {NewPolyBezierDerivEval-KeepDegree}{$n, t, r, \p{W}$}

\State Compute $h_0,h_1,\ldots,h_n$ as in \cite[Algorithm 2.2]{WCh2020}

\For {$k \gets 0,n$}
    \State $\p{u}^{(0)}_k \gets \p{W}_k$
\EndFor

\For {$j \gets 1,r$}
   \State $\p{u}^{(j)}_0 \gets 0$
    \For {$k \gets 0,n-1$}
        \State $\Delta \gets \p{u}^{(j-1)}_{k+1}-\p{u}^{(j-1)}_k$
        \State $\p{u}^{(j)}_k \gets \p{u}^{(j)}_k+(n-k)\cdot\Delta$
        \State $\p{u}^{(j)}_{k+1} \gets (k+1)\cdot\Delta$
    \EndFor
\EndFor

\For {$j \gets 0,r$}
    \State $\p{P}^{(j)}_0 \gets \p{u}^{(j)}_0$
    \For {$k \gets 1,n$}
        \State $\p{P}^{(j)}_k \gets (1-h_k) \cdot \p{P}^{(j)}_{k-1}
                                  + h_k \cdot \p{u}^{(j)}_k$  
    \EndFor
\EndFor

\State \Return $\p{P}^{(0)}_n,\p{P}^{(1)}_n,\ldots,\p{P}^{(n)}_n$

\EndProcedure
\end{algorithmic}
\end{algorithm}

It is worth mentioning again that when solving the polynomial version of 
Problem~\ref{P:ComputeDeriv} using Algorithm~\ref{A:ModPolyBezierDeriv}, one 
has to compute just $n+1$ quantities $h^n_0(t), h^n_1(t),\ldots, h^n_n(t)$
(cost $O(n)$) while in Algorithm~\ref{A:PolyBezierDeriv} the computation of 
$(r+1)(2n+2-r)/2$ quantities $h^{n-j}_i(t)$ for all $j=0,1,\ldots,r$ and all 
$i=0,1,\ldots,n-j$ is required (with $O(n^2)$ complexity if $r$ is $O(n)$) --- 
to do this, it is neccessary to call the procedure \textsc{NewPolyBezierEval} 
$r+1$ times.

%%%%%%%%%%%%%%%%%%%%%%%%%%%%%%%%%%%%%%%%%%%%%%%%%%%%%%%%%%%%%%%%%%%%%%%%%%%%%%
\subsection{Numerical tests}                          \label{SS:PolyDiff-Test}
%%%%%%%%%%%%%%%%%%%%%%%%%%%%%%%%%%%%%%%%%%%%%%%%%%%%%%%%%%%%%%%%%%%%%%%%%%%%%%

The results presented in Tables~\ref{T:Table1}--\ref{T:Table3} have been 
obtained on a computer with \texttt{Intel Core i5-6300U CPU} at \texttt{2.40GHz}
processor and \texttt{8GB} \texttt{RAM}, using \texttt{G++ 11.3.0} (double 
precision). The source code in \texttt{C++17} which was used to perform all 
tests is available at 
\url{http://www.ii.uni.wroc.pl/~pwo/programs/NewBezierDiffEval.cpp}.

\begin{example}\label{E:Example1}
Tables~\ref{T:Table1}, \ref{T:Table1-3D} and~\ref{T:Table2} show the comparison
between the running times of the de Casteljau algorithm and the new methods for
the evaluation of the derivatives of a B\'{e}zier curve proposed 
in~\S\ref{SS:PolyDiff-I} and~\S\ref{SS:PolyDiff-II}.  

The following numerical experiments have been conducted. For fixed $n$ and $r$ 
such that $n \geq r\geq 0$, $1000$ test sets consisting of one polynomial 
B\'{e}zier curve of degree $n$ are generated. Their control points 
$\p{W}_k\in[-1,1]^d$ $(0\leq k\leq n;\; d\in\{1,2,3\})$ have been generated 
using the \texttt{std::uniform\_real\_distribution} \texttt{C++17} function. The
value and the first $r$ derivatives of each curve are then evaluated at 501
points $t_i:=i/500$ $(0\leq i\leq 500)$. Each algorithm is tested using the same 
curves. Tables~\ref{T:Table1} $(d=2)$, \ref{T:Table1-3D} $(d=3)$
and~\ref{T:Table2} $(d=1)$ show the total running time of all $501\times 1000$
evaluations of the values and first $r$ derivatives.
\end{example}

\begin{table*}[ht!]
{\small
\begin{center}
\renewcommand{\arraystretch}{1.45}
\begin{tabular}{lcccc}
$n$ & $r$ & de~Casteljau & new method & new method (kept degree)\\ \hline
2 & 1 & 0.101701 & \textbf{0.085681} & 0.094308 \\
2 & 2 & 0.135102 & \textbf{0.096371} & 0.135967 \\ \hline
3 & 1 & 0.15811 & \textbf{0.131103} & 0.136471 \\
3 & 2 & 0.186423 & \textbf{0.161211} & 0.203425 \\
3 & 3 & 0.232719 & \textbf{0.170037} & 0.263575 \\ \hline
4 & 1 & 0.239295 & 0.178819 & \textbf{0.177125} \\
4 & 2 & 0.267154 & \textbf{0.226221} & 0.26364 \\
4 & 3 & 0.316511 & \textbf{0.255086} & 0.348704 \\
4 & 4 & 0.383342 & \textbf{0.268832} & 0.421284 \\ \hline
5 & 1 & 0.342637 & 0.221575 & \textbf{0.217411} \\
5 & 2 & 0.370283 & \textbf{0.297803} & 0.320969 \\
5 & 3 & 0.424055 & \textbf{0.351203} & 0.417808 \\
5 & 5 & 0.592397 & \textbf{0.388229} & 0.611776 \\ \hline
10 & 1 & 1.20014 & 0.475849 & \textbf{0.442431} \\
10 & 2 & 1.21925 & 0.665185 & \textbf{0.639336} \\
10 & 3 & 1.27028 & 0.848587 & \textbf{0.837634} \\
10 & 10 & 2.6407 & \textbf{1.37989} & 2.23118 \\ \hline
20 & 1 & 4.19273 & 0.933546 & \textbf{0.856092} \\
20 & 2 & 4.2114 & 1.35745 & \textbf{1.22954} \\
20 & 3 & 4.25464 & 1.76631 & \textbf{1.61163} \\
20 & 20 & 14.2004 & \textbf{5.07646} & 8.16745 \\ \hline
30 & 1 & 9.17692 & 1.4094 & \textbf{1.24239} \\
30 & 2 & 9.19503 & 2.0426 & \textbf{1.78735} \\
30 & 3 & 9.23536 & 2.70152 & \textbf{2.40574} \\
30 & 30 & 40.6902 & \textbf{11.1665} & 17.4894 \\ \hline
50 & 1 & 25.7738 & 2.32732 & \textbf{2.06813} \\
50 & 2 & 25.7696 & 3.45315 & \textbf{2.96614} \\
50 & 3 & 25.8017 & 4.56996 & \textbf{3.87413} \\
50 & 50 & 162.676 & \textbf{30.4085} & 46.3989 \\ \hline
300 & 1 & 1005.76 & 14.9495 & \textbf{13.1013} \\
300 & 2 & 980.503 & 22.4398 & \textbf{18.8115} \\
300 & 3 & 984.143 & 30.5473 & \textbf{24.5207} \\
300 & 300 & \mbox{max.~time exceeded} & \textbf{1108.95} & 1704.65\\
\end{tabular}
\renewcommand{\arraystretch}{1}
\vspace{2ex}
\caption{Running times comparison (in seconds) for Example~\ref{E:Example1} with
$d=2$.}\label{T:Table1}
\vspace{-3ex}
\end{center}
}
\end{table*}

\begin{table*}[ht!]
{\small
\begin{center}
\renewcommand{\arraystretch}{1.45}
\begin{tabular}{lcccc}
$n$ & $r$ & de~Casteljau & new method & new method (kept degree)\\ \hline
2 & 1 & 0.114709 & \textbf{0.077897} & 0.083363 \\
2 & 2 & 0.139229 & \textbf{0.0866} & 0.117128 \\ \hline
3 & 1 & 0.187191 & \textbf{0.117957} & 0.122946 \\
3 & 2 & 0.222727 & \textbf{0.139821} & 0.166657 \\
3 & 3 & 0.304963 & \textbf{0.147241} & 0.229243 \\ \hline
4 & 1 & 0.244661 & \textbf{0.160038} & 0.18525 \\
4 & 2 & 0.273558 & \textbf{0.204664} & 0.241133 \\
4 & 3 & 0.333152 & \textbf{0.219738} & 0.290539 \\
4 & 4 & 0.446606 & \textbf{0.261004} & 0.370166 \\ \hline
5 & 1 & 0.346794 & \textbf{0.189484} & 0.195475 \\
5 & 2 & 0.37635 & \textbf{0.278451} & 0.278773 \\
5 & 3 & 0.421566 & \textbf{0.29843} & 0.348953 \\
5 & 5 & 0.626254 & \textbf{0.330816} & 0.497813 \\ \hline
10 & 1 & 1.07925 & 0.399504 & \textbf{0.383105} \\
10 & 2 & 1.10854 & 0.563189 & \textbf{0.533415} \\
10 & 3 & 1.15732 & 0.704217 & \textbf{0.683903} \\
10 & 10 & 2.85077 & \textbf{1.16324} & 1.73712 \\ \hline
20 & 1 & 3.79591 & 0.793563 & \textbf{0.716656} \\
20 & 2 & 3.83577 & 1.15561 & \textbf{1.00517} \\
20 & 3 & 3.87833 & 1.49507 & \textbf{1.29573} \\
20 & 20 & 17.0293 & \textbf{4.27409} & 6.21786 \\ \hline
30 & 1 & 8.17304 & 1.18627 & \textbf{1.05438} \\
30 & 2 & 8.21327 & 1.74525 & \textbf{1.48712} \\
30 & 3 & 8.25839 & 2.28333 & \textbf{1.91689} \\
30 & 30 & 50.9817 & \textbf{9.39334} & 13.5276 \\ \hline
50 & 1 & 23.4891 & 1.99188 & \textbf{1.72466} \\
50 & 2 & 23.553 & 2.95276 & \textbf{2.43912} \\
50 & 3 & 23.648 & 3.86835 & \textbf{3.14137} \\
50 & 50 & 208.454 & \textbf{25.671} & 36.8387 \\ \hline
300 & 1 & 843.983 & 12.0198 & \textbf{10.2502} \\
300 & 2 & 844.239 & 17.9918 & \textbf{14.4725} \\
300 & 3 & 843.657 & 23.9849 & \textbf{18.7152} \\
300 & 300 & \mbox{max.~time exceeded} & \textbf{1128.07} & 1544.71\\
\end{tabular}
\renewcommand{\arraystretch}{1}
\vspace{2ex}
\caption{Running times comparison (in seconds) for Example~\ref{E:Example1} with
$d=3$.}\label{T:Table1-3D}
\vspace{-3ex}
\end{center}
}
\end{table*}

\begin{table*}[ht!]
{\small
\begin{center}
\renewcommand{\arraystretch}{1.45}
\begin{tabular}{lcccc}
$n$ & $r$ & de~Casteljau & new method & new method (kept degree)\\ \hline
20 & 1 & 0.424387 & 0.306281 & \textbf{0.205171} \\
20 & 2 & 0.433771 & 0.440159 & \textbf{0.24861} \\
20 & 3 & 0.450863 & 0.567658 & \textbf{0.297138} \\
20 & 20 & 3.6252 & 1.59332 & \textbf{1.01604} \\ \hline
30 & 1 & 0.915564 & 0.465021 & \textbf{0.307133} \\
30 & 2 & 0.925343 & 0.68182 & \textbf{0.371499} \\
30 & 3 & 0.942442 & 0.88906 & \textbf{0.439824} \\
30 & 30 & 11.4358 & 3.54949 & \textbf{2.1813} \\ \hline
50 & 1 & 2.47771 & 0.790672 & \textbf{0.514966} \\
50 & 2 & 2.51905 & 1.17153 & \textbf{0.630606} \\
50 & 3 & 2.5071 & 1.5438 & \textbf{0.740185} \\
50 & 50 & 48.8375 & 10.013 & \textbf{5.92621} \\ \hline
200 & 1 & 37.7228 & 3.24994 & \textbf{2.07921} \\
200 & 2 & 37.6743 & 4.84807 & \textbf{2.52765} \\
200 & 3 & 37.7775 & 6.46446 & \textbf{2.97356} \\
200 & 200 & 2548.32 & 162.713 & \textbf{91.0993} \\ \hline
300 & 1 & 83.895 & 4.95148 & \textbf{3.17461} \\
300 & 2 & 84.0311 & 7.42649 & \textbf{3.86009} \\
300 & 3 & 84.3917 & 9.85612 & \textbf{4.52994} \\
300 & 300 & 8498.56 & 368.917 & \textbf{203.697} \\
\end{tabular}
\renewcommand{\arraystretch}{1}
\vspace{2ex}
\caption{Running times comparison (in seconds) for Example~\ref{E:Example1} 
with $d=1$.}\label{T:Table2}
\vspace{-3ex}
\end{center}
}
\end{table*}

\begin{example}\label{E:Example1Many}
Table~\ref{T:Table3} shows the comparison between the running times of the de 
Casteljau algorithm and the new methods proposed in~\S\ref{SS:PolyDiff-I} 
and~\S\ref{SS:PolyDiff-II} when we evaluate the derivatives of multiple 
polynomial B\'{e}zier curves for the same values of a parameter $t\in[0,1]$ 
(cf.~Remark~\ref{R:Remark2}).  

The following numerical experiments have been conducted. For fixed $n$ and $r$ 
such that $n \geq r\geq 0$, $1000$ test sets consisting of $m$ polynomial 
B\'{e}zier curves of degree $n$ are generated. Their control points 
$\p{W}_k\in[-1,1]^2$ $(0\leq k\leq n)$ have been generated using the 
\texttt{std::uniform\_real\_distribution} \texttt{C++17} function. The value and
the first $r$ derivatives of each curve are then evaluated at 501 points 
$t_i:=i/500$ $(0\leq i\leq 500)$. Each algorithm is tested using the same 
curves. Table~\ref{T:Table3} shows the total running time of all 
$501\times 1000\times m$ evaluations of the values and first $r$
derivatives.
\end{example}

\begin{table*}[ht!]
{\small
\begin{center}
\renewcommand{\arraystretch}{1.45}
\begin{tabular}{lccccc}

$n$ & $r$ & $m$ & de~Casteljau & new method & new method (kept degree)\\ \hline
2 & 2 & 2 & 0.27069 & \textbf{0.166395} & 0.270576\\
  &   & 5 & 0.687538 & \textbf{0.382816} & 0.660179\\
  &   & 10 & 1.325 & \textbf{0.736941} & 1.27104\\ \hline

3 & 2 & 2 & 0.395928 & \textbf{0.289248} & 0.397414\\
  &   & 5 & 0.936259 & \textbf{0.669781} & 0.961091\\
  &   & 10 & 1.86281 & \textbf{1.29078} & 1.9029\\ \hline

3 & 3 & 2 & 0.486215 & \textbf{0.302551} & 0.515761\\
  &   & 5 & 1.1644 & \textbf{0.697202} & 1.25619\\
  &   & 10 & 2.31923 & \textbf{1.34652} & 2.50325\\ \hline

5 & 2 & 2 & 0.748109 & \textbf{0.545108} & 0.633346\\
  &   & 5 & 1.85355 & \textbf{1.25085} & 1.52479\\
  &   & 10 & 3.69939 & \textbf{2.47735} & 3.03078\\ \hline

5 & 5 & 2 & 1.1887 & \textbf{0.706878} & 1.22357\\
  &   & 5 & 2.95723 & \textbf{1.62799} & 2.98678\\
  &   & 10 & 5.91091 & \textbf{3.18372} & 5.92845\\ \hline

10 & 2 & 2 & 2.48081 & \textbf{1.18779} & 1.23888\\
   &   & 5 & 6.099 & \textbf{2.76847} & 3.02924\\
   &   & 10 & 12.1945 & \textbf{5.4478} & 6.04751\\ \hline

10 & 10 & 2 & 5.28712 & \textbf{2.4833} & 4.43912\\
   &    & 5 & 13.2077 & \textbf{5.84922} & 11.0316\\
   &    & 10 & 26.4987 & \textbf{11.4373} & 21.8572\\ \hline

50 & 2 & 2 & 51.6284 & \textbf{6.10374} & 5.67993\\
   &   & 5 & 128.862 & 14.144 & \textbf{13.9533}\\
   &   & 10 & 257.731 & \textbf{27.5226} & 27.7681\\ \hline

50 & 50 & 2 & 324.705 & \textbf{53.7153} & 92.7871\\
   &    & 5 & 813.543 & \textbf{124.272} & 233.833\\
   &    & 10 & 1624.88 & \textbf{243.215} & 466.352\\
\end{tabular}
\renewcommand{\arraystretch}{1}
\vspace{2ex}
\caption{Running times comparison (in seconds) for 
Example~\ref{E:Example1Many}, for $d=2$.}\label{T:Table3}
\vspace{-3ex}
\end{center}
}
\end{table*}

From the conducted tests, it follows that the numerical performance of the 
methods presented in Sections~\ref{SS:PolyDiff-I} and~\ref{SS:PolyDiff-II} was 
very close to that of the de Casteljau algorithm.

As expected, keeping the degree of a B\'{e}zier curve when finding a derivative
(cf.~\S\ref{SS:PolyDiff-II}) is slower than using the new method 
with the degrees of derivative vector curves kept as low as possible 
(cf.~\S\ref{SS:PolyDiff-I}) unless $n$ is quite big and $r$ is small. In both 
variants, the new method outperformed the de~Casteljau algorithm. See 
Tables~\ref{T:Table1}, \ref{T:Table1-3D} and~\ref{T:Table3}.

However, in Example~\ref{E:Example1}, if $d=1$ and $n$ is quite big the method 
from~\S\ref{SS:PolyDiff-II} is the fastest one. Again, both new methods 
outperformed the de~Casteljau algorithm. See Table~\ref{T:Table2}. 

%%%%%%%%%%%%%%%%%%%%%%%%%%%%%%%%%%%%%%%%%%%%%%%%%%%%%%%%%%%%%%%%%%%%%%%%%%%%%%
%%%%%%%%%%%%%%%%%%%%%%%%%%%%%%%%%%%%%%%%%%%%%%%%%%%%%%%%%%%%%%%%%%%%%%%%%%%%%%
\section{Derivatives of rational B\'{e}zier curves}          \label{S:RatDiff}
%%%%%%%%%%%%%%%%%%%%%%%%%%%%%%%%%%%%%%%%%%%%%%%%%%%%%%%%%%%%%%%%%%%%%%%%%%%%%%
%%%%%%%%%%%%%%%%%%%%%%%%%%%%%%%%%%%%%%%%%%%%%%%%%%%%%%%%%%%%%%%%%%%%%%%%%%%%%%

Using the well-known differentiation rules, the following general formula 
for derivatives of a rational B\'{e}zier curve~\eqref{E:RationalBezierCurve} can 
be obtained:
\begin{equation}\label{E:RatBezDiff}
\p{R}_n^{(k)}(t)=\dfrac{1}{A_0(t)}
        \left(
              \dfrac{d^k}{dt^k}\sum_{j=0}^n\omega_j\p{W}_jB^n_j(t)
                       -\sum_{i=0}^{k-1}\binom{k}{i}\p{R}_n^{(i)}(t)A_{k-i}(t)
        \right)\qquad (0\leq t\leq 1),
\end{equation}
where $k\in\mathbb N$, and
\begin{equation}\label{E:Def_A_i}
A_i(t):=\dfrac{d^i}{dt^i}\sum_{j=0}^n\omega_j B^n_j(t)\qquad (i=0,1,\ldots)
\end{equation}
(cf., e.g., \cite[Eq.~(13.9)]{Farin2002}).

In applications related to CAGD, it is usually sufficient to find the first,
second or third derivative of a parametric curve, e.g., to compute its curvature
and torsion or to merge curves smoothly. However, in this section, we consider 
Problem~\ref{P:ComputeDeriv} in its general setting remembering that the 
differentiation of rational functions is often necessary (for example,
in numerical analysis) and that the curve~\eqref{E:RationalBezierCurve} is a 
real-valued function if control points are in ${\mathbb E}^1$.

As in the polynomial case, we show that the use of the geometric 
scheme~\eqref{E:Def_h_k_Q_k} allows us to solve Problem~\ref{P:ComputeDeriv} 
very efficiently.

%%%%%%%%%%%%%%%%%%%%%%%%%%%%%%%%%%%%%%%%%%%%%%%%%%%%%%%%%%%%%%%%%%%%%%%%%%%%%%
\subsection{First method: accelerating Floater's formulas 
            for the first and the second derivative}
                                                    \label{SS:RatDiff-Floater}
%%%%%%%%%%%%%%%%%%%%%%%%%%%%%%%%%%%%%%%%%%%%%%%%%%%%%%%%%%%%%%%%%%%%%%%%%%%%%%

In~\cite{Floater91}, Floater gives useful formulas for the evaluation of the 
first two derivatives of rational B\'{e}zier curves which have a geometric 
interpretation and are closely related to the rational de Casteljau algorithm.

\begin{theorem}[{\cite[Proposition 2, Proposition 11]{Floater91}}]\label{T:ThmF}
Let $\p{R}_n : [0, 1]\rightarrow\mathbb{E}^d$ be a rational B\'{e}zier curve of
degree $n$ defined in~\eqref{E:RationalBezierCurve}. Then
\begin{equation}\label{E:FloaterDiff-1}
\p{R}_n^{(1)}(t)=n\dfrac{\omega^{(n-1)}_0\omega^{(n-1)}_1}
                        {\left(\omega^{(n)}_0\right)^2}\cdot 
                            \left(\p{W}^{(n-1)}_1-\p{W}^{(n-1)}_0\right),
\end{equation}
and
\begin{multline}\label{E:FloaterDiff-2}
\p{R}_n^{(2)}(t)=n\dfrac{\omega^{(n-2)}_2}{\left(\omega^{(n)}_0\right)^3}
                   \left(2n(\omega^{(n-1)}_0)^2-(n-1)\omega^{(n-2)}_0
                           \omega^{(n)}_0-2\omega^{(n-1)}_0\omega^{(n)}_0\right)
                   \left(\p{W}^{(n-2)}_2-\p{W}^{(n-2)}_1\right)\\                   
                -n\dfrac{\omega^{(n-2)}_0}{\left(\omega^{(n)}_0\right)^3}
                   \left(2n(\omega^{(n-1)}_1)^2-(n-1)\omega^{(n-2)}_2
                           \omega^{(n)}_0-2\omega^{(n-1)}_1\omega^{(n)}_0\right)
                   \left(\p{W}^{(n-2)}_1-\p{W}^{(n-2)}_0\right),
\end{multline}
where the numbers $\omega^{(i)}_k$ and points $\p{W}^{(i)}_k$ are computed by 
the rational de Casteljau algorithm (cf.~\eqref{E:RatdeCastel}).
\end{theorem}

\begin{figure}[ht!]
\centering
\begin{tikzpicture}
\node at (0, 5) {$\omega_0$};
\node at (0, 4) {$\omega_1$};
\node at (0, 3) {$\omega_2$};
\node at (0, 2) {$\omega_3$};
\node at (0, 1) {$\omega_4$};
\node at (0, 0) {$\omega_5$};

\node at (1.5, 5) {$\omega^{(1)}_0$};
\node at (1.5, 4) {$\omega^{(1)}_1$};
\node at (1.5, 3) {$\omega^{(1)}_2$};
\node at (1.5, 2) {$\omega^{(1)}_3$};
\node at (1.5, 1) {$\omega^{(1)}_4$};

\node at (3, 5) {$\omega^{(2)}_0$};
\node at (3, 4) {$\omega^{(2)}_1$};
\node at (3, 3) {$\omega^{(2)}_2$};
\node at (3, 2) {$\omega^{(2)}_3$};

\node at (4.5, 5) {$\omega^{(3)}_0$};
\node at (4.5, 4) {$\omega^{(3)}_1$};
\node at (4.5, 3) {$\omega^{(3)}_2$};

\node at (6, 5) {$\omega^{(4)}_0$};
\node at (6, 4) {$\omega^{(4)}_1$};

\node at (7.5, 5) {$\omega^{(5)}_0$};

\foreach \y in {0,...,4}
{
  \draw[->] (0.4, 5 - \y) -- (1.1, 5 - \y);
  \draw[->] (0.4, 4.2 - \y) -- (1.1, 4.8 - \y);
}

\foreach \y in {0,...,3}
{
  \draw[->] (1.9, 5 - \y) -- (2.6, 5 - \y);
  \draw[->] (1.9, 4.2 - \y) -- (2.6, 4.8 - \y);
}

\foreach \y in {0,...,2}
{
  \draw[->] (3.4, 5 - \y) -- (4.1, 5 - \y);
  \draw[->] (3.4, 4.2 - \y) -- (4.1, 4.8 - \y);
}

\foreach \y in {0,...,1}
{
  \draw[->] (4.9, 5 - \y) -- (5.6, 5 - \y);
  \draw[->] (4.9, 4.2 - \y) -- (5.6, 4.8 - \y);
}

\foreach \y in {0,...,0}
{
  \draw[->] (6.4, 5 - \y) -- (7.1, 5 - \y);
  \draw[->] (6.4, 4.2 - \y) -- (7.1, 4.8 - \y);
}
\end{tikzpicture}
\caption{The evaluation of the intermediate weights necessary for the base 
Floater's algorithm $(n=5)$. The arrows denote the steps of the rational de 
Casteljau algorithm. The same scheme applies for the necessary intermediate 
points $\p{W}^{(i)}_k$.}\label{F:Figure2}
\end{figure}
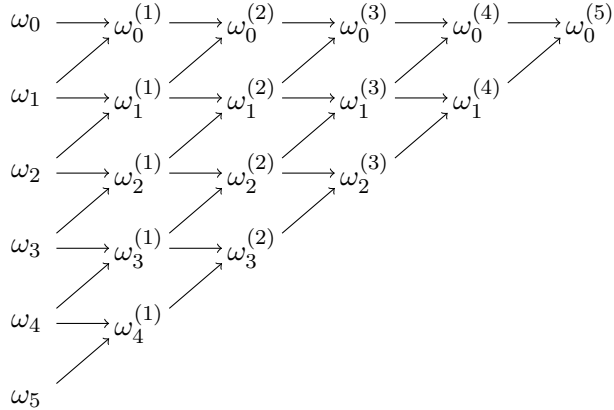

Classically, the vectors $\p{R}_n^{(1)}(t),\p{R}_n^{(2)}(t)\in{\mathbb R}^d$ 
given by~\eqref{E:FloaterDiff-1} and \eqref{E:FloaterDiff-2} could be found in 
$O(n^2d)$ time using the rational de Casteljau algorithm, which would also give 
the total computational complexity of finding the first two derivatives of 
a rational B\'{e}zier curve using the mentioned Floater's result (see also 
Figure~\ref{F:Figure2}).

Note that when evaluating the Floater's equations~\eqref{E:FloaterDiff-1} and 
\eqref{E:FloaterDiff-2}, one could avoid a significant part of the computations 
by using the linear-time geometric method proposed in~\cite{WCh2020} 
(see also the scheme~\eqref{E:Def_h_k_Q_k}). Indeed, the numbers 
$\omega^{(n-2)}_k$ and points $\p{W}^{(n-2)}_k$ ($k\in\{0,1,2\}$) are given by 
the formulas
$$
\omega^{(n-2)}_k=\sum_{i=0}^{n-2}\omega_{i+k}B^{n-2}_i(t),\qquad
\p{W}^{(n-2)}_k=\dfrac{\displaystyle
                         \sum_{i=0}^{n-2}\omega_{i+k}B^{n-2}_i(t)\p{W}_{i+k}}
                      {\displaystyle
                         \sum_{i=0}^{n-2}\omega_{i+k}B^{n-2}_i(t)}.
$$
See, e.g., \cite[\S13]{Farin2002}, \cite[Eq.~(3), (5)]{Floater91}.
This means that they can be treated as polynomial B\'{e}zier curves of degree 
$n-2$ with one-dimensional control points $\omega_{i+k}$ and rational B\'{e}zier
curves of degree $n-2$ with $d$-dimensional control points $\p{W}_{i+k}$, 
respectively. Instead of finding their values with the rational de Casteljau 
algorithm, using~\cite[{Algorithm~2.2}]{WCh2020} would give the necessary 
quantities with lower complexity. Afterwards, one can find the remaining 
necessary numbers $\omega^{(n-1)}_0,\omega^{(n-1)}_1,\omega^{(n)}_0$ and points 
$\p{W}^{(n-1)}_0,\p{W}^{(n-1)}_1$ using the recurrence relations from the 
rational de Casteljau algorithm. Certainly, this approach has a geometric 
interpretation.

The idea of the presented approach for acceleration of Floater's formulas can 
be seen in Figure~\ref{F:Figure3}. This gives the \textit{geometric} method for 
computing the first two derivatives $\p{R}_n^{(1)}(t)$ and 
$\p{R}_n^{(2)}(t)\in{\mathbb R}^d$ in $O(nd)$ time.

\begin{figure}[ht!]
\centering
\begin{tikzpicture}
\node at (0, 5) {$\omega_0$};
\node at (0, 4) {$\omega_1$};
\node at (0, 3) {$\omega_2$};
\node at (0, 2) {$\omega_3$};
\node at (0, 1) {$\omega_4$};
\node at (0, 0) {$\omega_5$};

\node at (4.5, 5) {$\omega^{(3)}_0$};
\node at (4.5, 4) {$\omega^{(3)}_1$};
\node at (4.5, 3) {$\omega^{(3)}_2$};

\node at (6, 5) {$\omega^{(4)}_0$};
\node at (6, 4) {$\omega^{(4)}_1$};

\node at (7.5, 5) {$\omega^{(5)}_0$};

\foreach \y in {0,...,1}
{
  \draw[->] (4.9, 5 - \y) -- (5.6, 5 - \y);
  \draw[->] (4.9, 4.2 - \y) -- (5.6, 4.8 - \y);
}

\foreach \y in {0,...,0}
{
  \draw[->] (6.4, 5 - \y) -- (7.1, 5 - \y);
  \draw[->] (6.4, 4.2 - \y) -- (7.1, 4.8 - \y);
}

\draw (-0.5, 5.4) -- (0.3, 5.4) -- (0.3, 1.6) -- (-0.5, 1.6) -- (-0.5, 5.4);
\draw (-0.4, 4.4) -- (0.4, 4.4) -- (0.4, 0.6) -- (-0.4, 0.6) -- (-0.4, 4.4);
\draw (-0.3, 3.4) -- (0.5, 3.4) -- (0.5, -0.4) -- (-0.3, -0.4) -- (-0.3, 3.4);

\draw[->>] (0.3, 5) -- (4.1, 5);
\draw[->>] (0.4, 4) -- (4.1, 4);
\draw[->>] (0.5, 3) -- (4.1, 3);

\end{tikzpicture}
\caption{The evaluation of the intermediate weights necessary for the 
accelerated Floater's algorithm $(n=5)$. The one-headed arrows denote the 
steps of the rational de Casteljau algorithm, while the two-headed arrows
represent the linear-time geometric method~\cite{WCh2020}. The same scheme 
applies for the necessary intermediate points $\p{W}^{(i)}_k$ 
(cf.~Figure~\ref{F:Figure2}).}\label{F:Figure3}
\end{figure}

%%%%%%%%%%%%%%%%%%%%%%%%%%%%%%%%%%%%%%%%%%%%%%%%%%%%%%%%%%%%%%%%%%%%%%%%%%%%%%
%
% Sprawdzić numery wzorów!!!
%
\subsection{\texorpdfstring{Second method: using 
                                        the scheme~\eqref{E:Def_h_k_Q_k} to 
                                        compute the derivative of any order}
                           {Second method: using the scheme (1.9) to 
                                        compute the derivative of any order}}  
                                                         \label{SS:RatDiff-II}
%%%%%%%%%%%%%%%%%%%%%%%%%%%%%%%%%%%%%%%%%%%%%%%%%%%%%%%%%%%%%%%%%%%%%%%%%%%%%%

Another possibility of solving Problem~\ref{P:ComputeDeriv} is the 
\textit{differentiation} of the geometric method from~\cite{WCh2020}. Using this
idea, one can derive differential-recurrence relations which allow us to find 
the $k$th derivative of rational B\'{e}zier curves for any natural
number $k$.

%%%%%%%%%%%%%%%%%%%%%%%%%%%%%%%%%%%%%%%%%%%%%%%%%%%%%%%%%%%%%%%%%%%%%%%%%%%%%%
\subsubsection{Derivatives of the first and the second order}   
                                                           \label{SS:F+S-Der}
%%%%%%%%%%%%%%%%%%%%%%%%%%%%%%%%%%%%%%%%%%%%%%%%%%%%%%%%%%%%%%%%%%%%%%%%%%%%%%

First, let us consider the problem of computing the first and the second
derivatives of the rational B\'{e}zier curve~\eqref{E:RationalBezierCurve} 
(cf.~\S\ref{SS:RatDiff-Floater}). 

Clearly, we have $\p{R}_n'(t)=(\p{Q}^n_n(t))'$ and 
$\p{R}_n''(t)=(\p{Q}^n_n(t))''$ (see~\eqref{E:Def_Q}). The differentiation of 
the relation for $\p{Q}_i\equiv\p{Q}^n_i(t)$ from the 
scheme~\eqref{E:Def_h_k_Q_k} gives
$$
\p{Q}_i'=\left((1-h_i)\p{Q}_{i-1}\right)'+h_i'\p{W}_i
             =(1-h_i)\p{Q}_{i-1}'+h_i'\left(\p{W}_i-\p{Q}_{i-1}\right)
                                                      \qquad (1\leq i\leq n).
$$
Here, $h_i\equiv h^n_i(t)$ (cf.~\eqref{E:BezierDefHSum}) and $\p{Q}_0'=\Theta$, 
where $\Theta$ is a zero vector. Repeating the process gives the relation
$$
\p{Q}_i''=h_i''\left(\p{W}_i-\p{Q}_{i-1}\right)
                       -2h_i'\p{Q}_{i-1}'+(1-h_i)\p{Q}_{i-1}''
                                                      \qquad (1\leq i\leq n),
$$
where $\p{Q}_0''=\Theta$. 

To evaluate these equations for a fixed $t\in[0,1]$, first, we have to compute 
numbers $h_i$ and points $\p{Q}_i$ using the scheme~\eqref{E:Def_h_k_Q_k}. 
Next, we have to compute $h_i'$ and $h_i''$. In order to find them, we will use
the following recurrence relation for $h_i$ which is equivalent to the one given
in Eq.~\eqref{E:Def_h_k_Q_k}:
\begin{equation}\label{E:RatBezhUseful}
h_i\left(\omega_{i-1}i(1-t)+\omega_it(n-i+1)h_{i-1}\right)=
                                              \omega_it(n-i+1)h_{i-1}.
\end{equation}
From this, we get, after the differentiation,
$$
h_i'\left(\omega_{i-1}i(1-t)+\omega_it(n-i+1)h_{i-1}\right)= 
      \omega_{i-1}ih_i+\omega_i(n-i+1)(1-h_i)\left(th_{i-1}'+h_{i-1}\right),
$$
or, equivalently,
$$
h_i'=\dfrac{h_i}{th_{i-1}}\left((1-h_i)(th_{i-1}'+h_{i-1})
                          +\dfrac{\omega_{i-1}}{\omega_i(n-i+1)}ih_i\right).
$$
Here $i=1,2,\ldots,n$, and $h_0'=0$. Let us also observe that the quantity 
$h_i/(th_{i-1})$, in fact, does not depend on $h_i$, because
\begin{equation}\label{E:Rel_h_i_h_i-1}
\dfrac{h_i}{th_{i-1}}=
        \dfrac{\omega_i(n-i+1)}
              {(1-t)\omega_{i-1}i+t(n-i+1)h_{i-1}}=:f^n_i(t)\equiv f_i
                                                     \qquad(1\leq i\leq n).
\end{equation}

One can get the relation for the second derivative in a similar manner, by 
starting from~\eqref{E:RatBezhUseful} and applying the second derivative to 
both sides. After some algebra, one gets
$$
h_i''=f_i\left((1-h_i)\left(th_{i-1}''+2h_{i-1}'\right)
                                  -2h_i'\left(th_{i-1}'+h_{i-1}\right)+
                                   \dfrac{2\omega_{i-1}}{\omega_i(n-i+1)}ih_i'
                            \right),
$$
where $1\leq i\leq n$, and $h_0''=0$.

Thus,
\begin{equation}\label{E:SchemeDiff-1}
  \left\{
    \begin{array}{l}
      h_0':=0, \qquad \p{Q}_0':=\Theta, \\[1ex]
      h_i':=f_i\left((1-h_i)(th_{i-1}'+h_{i-1})
                   +\dfrac{\omega_{i-1}}{\omega_i(n-i+1)}ih_i\right),\\[2.75ex]
      \p{Q}_i':=(1-h_i)\p{Q}_{i-1}'+h_i'\left(\p{W}_i-\p{Q}_{i-1}\right),
    \end{array}
  \right.
\end{equation}
with $i=1,2,\ldots,n$ and $\Theta$ being the zero vector, and
\begin{equation}\label{E:SchemeDiff-2}
  \left\{
    \begin{array}{l}
      h_0'':=0, \qquad \p{Q}_0'':=\Theta, \\[1ex]    
      h_i'':=f_i\left(
                              (1-h_i)\left(th_{i-1}''+2h_{i-1}'\right)
                              -2h_i'\left(th_{i-1}'+h_{i-1}\right)+
                                 \dfrac{2\omega_{i-1}}
                                       {\omega_i(n-i+1)}ih_i'\right),\\[2.75ex]
      \p{Q}_i'':=h_i''\left(\p{W}_i-\p{Q}_{i-1}\right)
                                  -2h_i'\p{Q}_{i-1}'+(1-h_i)\p{Q}_{i-1}'',
    \end{array}
  \right.
\end{equation}
where $i=1,2,\ldots,n$ (cf.~\eqref{E:Rel_h_i_h_i-1}). Then 
$\p{R}_n'(t)=\p{Q}_n'$ and $\p{R}_n''(t)=\p{Q}_n''$ $(0\leq t\leq 1)$.

Taking into account that numbers $h_i$ and points $\p{Q}_i\in{\mathbb E}^d$ for 
$i=0,1,\ldots,n$ can be computed by the scheme~\eqref{E:Def_h_k_Q_k} in 
a linear time, one can find the first and next the second derivative of 
a rational B\'{e}zier curve using recurrence relations~\eqref{E:SchemeDiff-1} 
and~\eqref{E:SchemeDiff-2} in $O(nd)$ time. Note that this method is 
\textit{geometric}, because one only computes convex combination of points 
(cf.~\eqref{E:Def_h_k_Q_k}) and linear combinations of vectors in each step. 

%%%%%%%%%%%%%%%%%%%%%%%%%%%%%%%%%%%%%%%%%%%%%%%%%%%%%%%%%%%%%%%%%%%%%%%%%%%%%%
\subsubsection{Higher order derivatives}          \label{SS:RatDiff-II-Higher}
%%%%%%%%%%%%%%%%%%%%%%%%%%%%%%%%%%%%%%%%%%%%%%%%%%%%%%%%%%%%%%%%%%%%%%%%%%%%%%

For a fixed $t\in[0,1]$, to find the higher order derivatives of the quantities 
$h_i\equiv h_i^n(t)$ and the points 
$\p{Q}_i\equiv \p{Q}^n_i(t)\in{\mathbb E}^d$, one can use the well-known formula
$$
\left(f(t)\cdot g(t)\right)^{(r)}=
       \sum_{j=0}^r\binom{r}{j}f^{(j)}(t)g^{(r-j)}(t)\qquad (r\in{\mathbb N}).
$$

More precisely, for $\p{Q}_i$, starting from Eq.~\eqref{E:Def_h_k_Q_k}, one gets
\begin{eqnarray*}
\p{Q}_i^{(r)}&=&\left((1-h_i)\p{Q}_{i-1}\right)^{(r)}+h_i^{(r)}\p{W}_i
              =\sum_{j=0}^r\binom{r}{j}(1-h_i)^{(j)}\p{Q}_{i-1}^{(r-j)}+
                                                          h_i^{(r)}\p{W}_i\\
             &=&h_i^{(r)}\left(\p{W}_i-\p{Q}_{i-1}\right)+\p{Q}_{i-1}^{(r)}
                    -\sum_{j=0}^{r-1}\binom{r}{j}h_i^{(j)}\p{Q}_{i-1}^{(r-j)}.
\end{eqnarray*}

For $h_i$, starting from Eq.~\eqref{E:RatBezhUseful} gives
$$
\sum_{j=0}^r\binom{r}{j}h_i^{(j)}\left(\omega_{i-1}i(1-t)+
                                    \omega_ih_{i-1}t(n-i+1)\right)^{(r-j)}= 
                \omega_i(n-i+1)\sum_{j=0}^r\binom{r}{j}t^{(j)}h_{i-1}^{(r-j)},
$$
or, equivalently,
$$
h_i^{(r)}=f_i
             \left(
                th_{i-1}^{(r)}+rh_{i-1}^{(r-1)}
                    +\dfrac{\omega_{i-1}i}{\omega_i(n-i+1)}rh_i^{(r-1)}
                               -\sum_{j=1}^{r}\binom{r}{j}h_i^{(r-j)}
                                   \left(th_{i-1}^{(j)}+jh_{i-1}^{(j-1)}\right) 
             \right)
$$
(see also~\eqref{E:Rel_h_i_h_i-1}).

Thus, the general differential-recurrence relation which---together with the 
scheme~\eqref{E:Def_h_k_Q_k}---allows us to compute derivatives of the rational 
B\'{e}zier curve~\eqref{E:RationalBezierCurve} consecutively is as follows:
\begin{equation*}
\left\{
  \begin{array}{l}
  \displaystyle
  h^{(k)}_0:=0,\qquad \p{Q}_0^{(k)}:=\Theta,\\[1ex]
  \displaystyle
  h_i^{(k)}:=f_i
             \left(
                   (1-h_i)g_{k,i-1}
                    +\dfrac{\omega_{i-1}i}{\omega_i(n-i+1)}kh_i^{(k-1)}
                           -\sum_{j=1}^{k-1}\binom{k}{j}h_i^{(k-j)}g_{j,i-1}
             \right),\\[2.75ex]
  \displaystyle             
  \p{Q}_i^{(k)}:=h_i^{(k)}\left(\p{W}_i-\p{Q}_{i-1}\right)+\p{Q}_{i-1}^{(k)}
                    -\sum_{j=0}^{k-1}\binom{k}{j}h_i^{(j)}\p{Q}_{i-1}^{(k-j)},             
  \end{array}
\right.
\end{equation*}
where $k=1,2,\ldots,r$, $i=1,2,\ldots,n$, and
$g_{j,i-1}:=th_{i-1}^{(j)}+jh_{i-1}^{(j-1)}$ $(1\leq j\leq k)$
(cf.~\eqref{E:Rel_h_i_h_i-1}). Then
$\p{R}_n^{(k)}(t)=\p{Q}_n^{(k)}$ $(k=1,2,\ldots,r)$. See 
also~\eqref{E:SchemeDiff-1} and~\eqref{E:SchemeDiff-2}. 

In order to solve Problem~\ref{P:ComputeDeriv}, it means to compute the point 
$\p{R}_n(t)\in{\mathbb E}^d$ and the vectors 
$\p{R}_n'(t),\p{R}_n''(t),\ldots,\p{R}_n^{(r)}(t)\in{\mathbb R}^d$ 
$(0\leq t\leq 1)$ using the presented approach, one needs to do $O(ndr^2)$ 
arithmetic operations.

%%%%%%%%%%%%%%%%%%%%%%%%%%%%%%%%%%%%%%%%%%%%%%%%%%%%%%%%%%%%%%%%%%%%%%%%%%%%%%
%
% Sprawdzić numery wzorów!!!
%
\subsection{\texorpdfstring{Third method: using Eq.~\eqref{E:RatBezDiff} 
                            to compute derivatives of any order}
                           {Third method: using Eq. (3.1) to compute 
                                                    derivatives of any order}}
                                                        \label{SS:RatDiff-III}
%%%%%%%%%%%%%%%%%%%%%%%%%%%%%%%%%%%%%%%%%%%%%%%%%%%%%%%%%%%%%%%%%%%%%%%%%%%%%%

In this section, we discuss the possibility of using Eq.~\eqref{E:RatBezDiff}
to find the $k$th derivative of the rational B\'{e}zier 
curve~\eqref{E:RationalBezierCurve} for the fixed $t\in[0,1]$ and 
$k\in\mathbb N$. 

Let us fix the natural number $k>0$. One can write the mentioned formula as 
follows:
\begin{equation}\label{E:Expresion-diff-R_n}
\p{R}_n^{(k)}(t)=\dfrac{1}{A_0(t)}
        \left(
              \sum_{j=0}^n\omega_j\p{W}_j\dfrac{d^k}{dt^k}B^n_j(t)-
                     A_k(t)\p{R}_n(t)
                       -\sum_{i=1}^{k-1}\binom{k}{i}\p{R}_n^{(i)}(t)A_{k-i}(t)
        \right),
\end{equation}
where $A_i(t)$ $(i=1,2,\ldots,k)$ is the $i$th derivative of the polynomial  
$A_0(t)=\sum_{j=0}^n\omega_j B^n_j(t)$ (cf.~\eqref{E:Def_A_i}).

Certainly, $A_0$ can be interpreted as a polynomial B\'{e}zier curve of 
degree $n$ with control points 
$\omega_0,\omega_1,\ldots,\omega_n\in\mathbb{E}^1$. Similarly, $A_i$ 
$(1\leq i\leq k\leq n)$ is a polynomial B\'{e}zier curve of degree $n-i$. One 
can evaluate all the necessary B\'{e}zier curves $A_i(t)$ 
$(i=0,1,\ldots,k\leq n)$ in total $O(kn)$ time. If $k\leq n$ this is exactly 
the problem considered in Sections~\ref{SS:PolyDiff-I} and \ref{SS:PolyDiff-II} 
for $d=1$ and $r=k$ (see also the results of the numerical test presented
in~\S\ref{SS:PolyDiff-Test}). If $k>n$ then $A_k(t)\equiv0$. Note that the 
computed values $A_i(t)$ can be re-used for different $k$.

Note that for a given $t\in[0,1]$, one can evaluate 
all the needed numbers $\dfrac{d^k}{dt^k}B^n_j(t)$ for all $0\leq j\leq n$ and 
$k\leq n$ using~\eqref{E:BernsteinDeriv} in $O(kn)$ time.

Now, let us look into the geometricity of 
expression~\eqref{E:Expresion-diff-R_n}. It is clear that for $k>1$
$$
\sum_{i=1}^{k-1}\binom{k}{i}\p{R}^{(i)}_n(t)A_{k-i}(t)
$$
is the linear combination of vectors $\p{R}^{(i)}_n(t)\in{\mathbb R}^d$ 
$(i=1,2,\ldots,k-1)$ and thus is itself a vector in ${\mathbb R}^d$. 

If $A_k(t)=0$ $(k>0)$, the quantity $A_k(t)\p{R}_n(t)$ is a zero vector. 
Moreover, for $k\leq n$, the sum 
$$
\sum_{j=0}^n\omega_j\p{W}_j\dfrac{d^k}{dt^k}B^n_j(t)
$$
is a vector from ${\mathbb R}^d$, because it is a linear combination of points 
from ${\mathbb E}^d$ with coefficients which sum up to $0=A_k(t)$. More 
precisely, we have
\begin{equation}\label{E:Def_V_k}
\sum_{j=0}^n\omega_j\p{W}_j\dfrac{d^k}{dt^k}B^n_j(t)=
  \sum_{j=1}^n\omega_j\dfrac{d^k}{dt^k}B^n_j(t)\left(\p{W}_j-\p{W}_0\right)
                                            =:\p{V}_k(t)\qquad (1\leq k\leq n).
\end{equation}

Thus, if $A_k(t)=0$ $(k>0)$, we have
$$
\p{R}_n^{(k)}(t)=\dfrac{1}{A_0(t)}
        \left(\p{V}_k(t)
                -\sum_{i=1}^{k-1}\binom{k}{i}\p{R}_n^{(i)}(t)A_{k-i}(t)\right)
$$
for $k\leq n$ and 
\begin{equation}\label{E:Diff_R_n_for_k>n}
\p{R}_n^{(k)}(t)=-\dfrac{1}{A_0(t)}
                     \sum_{i=k-n}^{k-1}\binom{k}{i}\p{R}_n^{(i)}(t)A_{k-i}(t)
\end{equation}
for $k>n$.

If $A_k(t)\neq0$ it implies that $1\leq k\leq n$ and we have
$$
\p{R}_n^{(k)}(t)=\dfrac{A_k(t)}{A_0(t)}\left(\p{D}_k(t)-\p{R}_n(t)\right)
   -\dfrac{1}{A_0(t)}\sum_{i=1}^{k-1}\binom{k}{i}\p{R}_n^{(i)}(t)A_{k-i}(t),
$$        
where the point $\p{D}_k(t)\in{\mathbb E^d}$ is defined by
$$
\p{D}_k(t):=\dfrac{\sum_{j=0}^n\omega_j\p{W}_j\dfrac{d^k}{dt^k}B^n_j(t)}
                  {\sum_{j=0}^n\omega_j\dfrac{d^{k}}{dt^{k}}B^n_j(t)}=
            \dfrac{1}{A_k(t)}
                   \sum_{j=0}^n\omega_j\p{W}_j\dfrac{d^k}{dt^k}B^n_j(t)
                                                        \qquad (1\leq k\leq n).
$$
Moreover, it could be converted to a rational B\'{e}zier curve of degree $n-k$, 
however this would require us to compute $n-k$ linear combinations of $k$ 
$d$-dimensional control points, which would be a significant burden. Instead, 
it would be cheaper to simply compute the necessary derivatives
of Bernstein polynomials and then evaluate the point $\p{D}_k(t)$. In the last 
step, it is worth to use the geometric method for the evaluation of 
\textit{rational parametric objects} of the form~\eqref{E:ParObject} associated 
with the \textit{basis functions} satisfying the conditions~\eqref{E:BasisF} 
(for details, see~\cite[\S1]{WCh2020}).   

In order to use~\cite[{Algorithm~1.1}]{WCh2020} properly, let us introduce the 
functions
\begin{equation}\label{E:Def_b^nk_j}
b^{(n,k)}_j(t):=\omega_j\dfrac{d^{k}}{dt^{k}}B^n_j(t),
\end{equation}
where $j=0,1,\ldots,n$ and $1\leq k\leq n$. Recall that for a given $t\in[0,1]$, 
one can evaluate all these functions using~\eqref{E:BernsteinDeriv} in $O(kn)$ 
time.  

Let us fix $t\in[0,1]$ and let us define
$$
B^{+}_k(t):=\sum_{\substack{j=0\\ b^{(n,k)}_j(t) \geq 0}}^n b^{(n,k)}_j(t),
\qquad
B^{-}_k(t):=\sum_{\substack{j=0\\ b^{(n,k)}_j(t) < 0}}^n
                                               \left|b^{(n,k)}_j(t)\right|
                                                       \qquad (1\leq k\leq n).                                               
$$
Certainly, $B^{+}_k(t)-B^{-}_k(t)=A_k(t)$ and it is impossible that 
$B^{+}_k(t)=B^{-}_k(t)=0$, since $A_k(t)\neq0$. 

If $B^{+}_k(t),B^{-}_k(t)\neq0$ then
\begin{equation}\label{E:Point_D_II}
\p{D}_k(t)=\dfrac{B^{+}_k(t)}{A_k(t)}\p{D}^{+}_k(t)-
           \dfrac{B^{-}_k(t)}{A_k(t)}\p{D}^{-}_k(t),
\end{equation}
where the points $\p{D}^{\pm}_k(t)\in{\mathbb E^d}$ $(1\leq k\leq n)$ are 
defined by
$$
\p{D}^{+}_k(t):=\dfrac{1}{B^{+}_k(t)}
       \sum_{\substack{j=0\\ b^{(n,k)}_j(t)\geq0}}^nb^{(n,k)}_j(t)\p{W}_j,
\qquad
\p{D}^{-}_k(t):=\dfrac{1}{B^{-}_k(t)}
  \sum_{\substack{j=0\\ b^{(n,k)}_j(t)<0}}^n\left|b^{(n,k)}_j(t)\right|\p{W}_j.
$$
If one of the numbers $B^{+}_k(t),B^{-}_k(t)$ is equal to zero, we simply omit 
the appropriate summand in~\eqref{E:Point_D_II}.

Both points $\p{D}^{\pm}_k(t)\in{\mathbb E^d}$ can be treated as 
\textit{rational parametric objects} of the form~\eqref{E:ParObject} associated 
with the basis functions
$$
b^{(n,k)}_j(t)\Big/B^{+}_k(t)
\qquad\mbox{or}\qquad
\left|b^{(n,k)}_j(t)\right|\Big/B^{-}_k(t)
$$
for all weights equal to $1$ (cf.~\eqref{E:BasisF}). 
One can evaluate these points using~\cite[{Algorithm~1.1}]{WCh2020} in $O(nd)$ 
time.

To sum up, in order to find consequtively the point 
$\p{R}_n(t)\in{\mathbb E^d}$ and the vectors 
$$
\p{R}^{(1)}_n(t),\p{R}^{(2)}_n(t),\ldots,\p{R}^{(r)}_n(t)\in{\mathbb R^d}
$$ 
for a given $t\in[0,1]$ and $r\in{\mathbb N}$, i.e., to solve 
Problem~\ref{P:ComputeDeriv} using Eq.~\eqref{E:RatBezDiff}, one has to perform 
at most $O(rd(n+r))$ arithmetic operations. 

For example, if $r\leq n$ (if $r>n$ see also~\eqref{E:Diff_R_n_for_k>n}) we need
to compute:
\begin{itemize}
\itemsep 1ex

\item the coefficients $\displaystyle \binom{k}{i}\in{\mathbb N}$ for 
      $k=0,1,\ldots,r$ and $i=0,1,\ldots,k$ by using Pascal's triangle --- total
      time $O(r^2)$;

\item the quantities $\displaystyle b^{(n,k)}_j(t)\in{\mathbb R}$ 
      (see~\eqref{E:Def_b^nk_j}) for $k=0,1,\ldots,r$ and $j=0,1,\ldots,n$ --- 
      total time $O(rn)$ (cf.~\eqref{E:BernsteinDeriv});

\item the value $A_0(t)$ and the derivatives $A_i(t)\in{\mathbb R}$ 
      (see~\eqref{E:Def_A_i}) for $i=1,\ldots,r$ --- total time $O(rn)$ 
      (cf.~Sections~\ref{SS:PolyDiff-I} and~\ref{SS:PolyDiff-II}; see
      also the results of the numerical test presented 
      in~\S\ref{SS:PolyDiff-Test});      
            
\item the point $\p{R}_n(t)\in{\mathbb E^d}$ (see~\eqref{E:RationalBezierCurve})
      using~\cite[{Algorithm~2.1}]{WCh2020} 
      (cf.~also the scheme~\eqref{E:Def_h_k_Q_k}) --- $O(nd)$ time;      
      
\item for all $0\leq k\leq r$ such that $A_k(t)\neq 0$, the points 
      $\p{D}_k(t)\in{\mathbb E^d}$ (cf.~\eqref{E:Point_D_II}) 
      using~\cite[{Algorithm~1.1}]{WCh2020} --- total time at most $O(rnd)$;
      
\item for all $0\leq k\leq r$ such that $A_k(t)=0$, the vectors 
      $\p{V}_k(t)\in{\mathbb R^d}$ (cf.~\eqref{E:Def_V_k}) --- total time 
      at most $O(rnd)$;

\item the vectors $\displaystyle \sum_{i=1}^{k-1}
                    \binom{k}{i}\p{R}_n^{(i)}(t)A_{k-i}(t)\in{\mathbb R^d}$ for 
      $k=2,3,\ldots,r$ --- total time $O(dr^2)$ assuming that all the needed 
      derivatives $\p{R}_n^{(i)}(t)$ have been just evaluated.
      
\end{itemize}

%%%%%%%%%%%%%%%%%%%%%%%%%%%%%%%%%%%%%%%%%%%%%%%%%%%%%%%%%%%%%%%%%%%%%%%%%%%%%%
\subsection{Numerical tests}                           \label{SS:RatDiff-Test}
%%%%%%%%%%%%%%%%%%%%%%%%%%%%%%%%%%%%%%%%%%%%%%%%%%%%%%%%%%%%%%%%%%%%%%%%%%%%%%

The results presented in this section have been obtained on a computer with
\texttt{Intel Core i5-6300U CPU} at \texttt{2.40GHz} processor and \texttt{8GB}
\texttt{RAM}, using \texttt{G++ 11.3.0} (double precision). The source code in
\texttt{C++17} which was used to perform all tests is available at 
\url{http://www.ii.uni.wroc.pl/~pwo/programs/NewBezierDiffEval.cpp}. 

For the first two derivatives of the rational B\'{e}zier curves, the results 
given by the original formulas derived by Floater (cf.~Theorem~\ref{T:ThmF}) 
were treated as a baseline. All three new methods 
(cf.~\S\ref{SS:RatDiff-Floater}, \S\ref{SS:F+S-Der} and 
\S\ref{SS:RatDiff-III} in the case of the first and the second 
derivative) have given results which were numerically very close 
to it.

For the higher derivatives ($r=3,4,\ldots$), due to a lack of a good enough 
baseline, more rigorous testing of the methods from 
Sections~\ref{SS:RatDiff-II-Higher} and~\ref{SS:RatDiff-III} was required. 
These methods were implemented also in \textsf{Maple{\small\texttrademark}~14} 
and tested in two different precisions (with \texttt{Digits:=16} and 
\texttt{Digits:=64}). This has shown that the method from 
\S\ref{SS:RatDiff-II-Higher} is losing its numerical stability when we 
compute higher order derivatives at $t\approx1$. The same was not the case for 
$t\approx0$, therefore a symmetry property has been used, i.e.,  the order of 
the control points and weights was reversed and the evaluation proceeded for 
$1-t$. After such modification, the numerical behavior of the two methods was 
similar, with the one using the derivative of the scheme~\eqref{E:Def_h_k_Q_k} 
having a higher mean number of preserved digits but at a cost of a lower 
minimum.

\begin{example}\label{E:Example2}
Table~\ref{T:Table4} shows the comparison between the running times of
the basic and accelerated Floater method for finding the first two derivatives
of rational B\'{e}zier curves (cf.~\S\ref{SS:RatDiff-Floater}), as well as the 
methods from \S\ref{SS:F+S-Der} and Section~\ref{SS:RatDiff-III}, for the first
and the second derivative.

The following numerical experiments have been conducted. For fixed $n$ and 
$r\in\{1, 2\}$, $1000$ curves of degree $n$ are generated. Their control points 
$\p{W}_k\in[-1,1]^2$ and weights $\omega_k\in[0.01,2]$ $(0\leq k\leq n)$ have 
been generated using the \texttt{std::uniform\_real\_distribution} 
\texttt{C++17} function. The value and the first $r$ derivatives of each curve 
are then evaluated at 501 points $t_i:=i/500$ $(0\leq i\leq 500)$. Each method 
is tested using the same curves. Table~\ref{T:Table3} shows the total running 
time of all $501\times 1000$ evaluations of the values and first $r$
derivatives.
\end{example}

\begin{table*}[ht!]
{\small
\begin{center}
\renewcommand{\arraystretch}{1.45}
\begin{tabular}{lccccc}
$n$ & $r$ & Floater & accelerated Floater & method from~\S\ref{SS:RatDiff-III} 
& method from~\S\ref{SS:F+S-Der}\\ \hline
2 & 1 & 0.748014 & \textbf{0.739513} & 0.913747 & 0.797715\\
2 & 2 & 0.784439 & \textbf{0.777955} & 1.10084 & 0.944074\\ \hline
3 & 1 & 0.814712 & \textbf{0.786472} & 0.984628 & 0.864872\\
3 & 2 & \textbf{0.855245} & 0.856333 & 1.21017 & 1.08168\\ \hline
4 & 1 & 0.915191 & \textbf{0.842004} & 1.05823 & 0.945162\\
4 & 2 & 0.952852 & \textbf{0.931909} & 1.32864 & 1.22646\\ \hline
5 & 1 & 1.03542 & \textbf{0.896145} & 1.13006 & 1.01117\\
5 & 2 & 1.08462 & \textbf{1.01752} & 1.44696 & 1.36547\\ \hline
6 & 1 & 1.1814 & \textbf{0.94614} & 1.20178 & 1.07673\\
6 & 2 & 1.22655 & \textbf{1.0889} & 1.55127 & 1.49089\\ \hline
7 & 1 & 1.35605 & \textbf{0.99957} & 1.2813 & 1.15015\\
7 & 2 & 1.40184 & \textbf{1.16418} & 1.67509 & 1.63017\\ \hline
8 & 1 & 1.54516 & \textbf{1.06691} & 1.35524 & 1.21758\\
8 & 2 & 1.59299 & \textbf{1.24584} & 1.80184 & 1.75975\\ \hline
9 & 1 & 1.75172 & \textbf{1.11961} & 1.42808 & 1.28634\\
9 & 2 & 1.80111 & \textbf{1.34761} & 1.91147 & 1.89085\\ \hline
10 & 1 & 1.98636 & \textbf{1.17635} & 1.50722 & 1.35828\\
10 & 2 & 2.03354 & \textbf{1.42976} & 2.03142 & 2.01061\\ \hline
% 15 & 1 & 3.48201 & 1.44821 & 1.89256 & 1.69082\\
% 15 & 2 & 3.52349 & 1.83579 & 2.60412 & 2.65991\\ \hline
20 & 1 & 5.51608 & \textbf{1.71934} & 2.27297 & 2.01723\\
20 & 2 & 5.55209 & \textbf{2.24295} & 3.1866 & 3.29914\\ \hline
30 & 1 & 11.1122 & \textbf{2.25864} & 3.0148 & 2.67847\\
30 & 2 & 11.1391 & \textbf{3.05858} & 4.28363 & 4.57151\\ \hline
50 & 1 & 28.9998 & \textbf{3.35259} & 4.41841 & 3.99799\\
50 & 2 & 29.0193 & \textbf{4.71003} & 6.47314 & 7.12824\\ \hline
100 & 1 & 114.485 & \textbf{6.23003} & 7.98442 & 7.7865 \\
100 & 2 & 114.021 & \textbf{8.97354} & 11.7484 & 14.0026 \\ \hline
200 & 1 & 448.431 & \textbf{11.77} & 15.0267 & 14.8115 \\
200 & 2 & 448.572 & \textbf{17.3473} & 22.4004 & 27.5966 \\ \hline
300 & 1 & 1004.37 & \textbf{17.5306} & 22.3993 & 22.3232 \\
300 & 2 & 1004.89 & \textbf{26.0027} & 33.4608 & 41.7452 \\
\end{tabular}
\renewcommand{\arraystretch}{1}
\vspace{2ex}
\caption{Running times comparison (in seconds) for Example~\ref{E:Example2}
(the first two derivatives, $d=2$).}
\label{T:Table4}
\vspace{-3ex}
\end{center}
}
\end{table*}

Example~\ref{E:Example2} shows that the the new approach proposed in 
\S\ref{SS:RatDiff-Floater} is the fastest of the considered methods, except one 
case for $n=3$ and $r=2$, where the classic Floater method is faster by 0.1\%. 

\begin{example}\label{E:Example3}
Table~\ref{T:Table5} shows the comparison between the running times of
the methods from Section~\ref{SS:RatDiff-II-Higher} and
Section~\ref{SS:RatDiff-III} for computing higher order derivatives of 
a rational B\'{e}zier curve.

The following numerical experiments have been conducted. For fixed $n$ and $r$ 
such that $n\geq r\geq 0$, $1000$ curves of degree $n$ are generated. Their 
control points $\p{W}_k\in[-1,1]^2$ and weights $\omega_k\in[0.01,2]$ 
$(0\leq k\leq n)$ have been generated using the 
\texttt{std::uniform\_real\_distribution} \texttt{C++17} function. The value and
the first $r$ derivatives of each curve are then evaluated at 501 points 
$t_i:=i/500$ $(0\leq i\leq 500)$. Each method is tested using the same curves. 
Table~\ref{T:Table5} shows the total running time of all $501\times 1000$ 
evaluations of the values and first $r$ derivatives.
\end{example}

\begin{table*}[ht!]
{\small
\begin{center}
\renewcommand{\arraystretch}{1.45}
\begin{tabular}{lccc}
$n$ & $r$ & method from~\S\ref{SS:RatDiff-III} 
                                  & method from~\S\ref{SS:RatDiff-II}\\ \hline
3 & 3 & 1.43675 & \textbf{1.35709}\\ \hline
4 & 3 & 1.59043 & \textbf{1.57812}\\
4 & 4 & \textbf{1.86569} & 2.00114\\ \hline
5 & 3 & \textbf{1.73728} & 1.78502\\
5 & 4 & \textbf{2.05224} & 2.2999\\
5 & 5 & \textbf{2.37311} & 2.90818\\ \hline
% 6 & 3 & 1.88543 & 1.99019\\
% 6 & 4 & 2.23994 & 2.60773\\
% 6 & 5 & 2.59816 & 3.33085\\ \hline
% 7 & 3 & 2.04128 & 2.20985\\
% 7 & 4 & 2.43765 & 2.92096\\
% 7 & 5 & 2.82351 & 3.75532\\ \hline
% 8 & 3 & 2.20247 & 2.41517\\
% 8 & 4 & 2.62883 & 3.22843\\
% 8 & 5 & 3.06815 & 4.16938\\ \hline
% 9 & 3 & 2.35178 & 2.62688\\
% 9 & 4 & 2.81991 & 3.53352\\
% 9 & 5 & 3.29853 & 4.59529\\ \hline
10 & 3 & \textbf{2.52806} & 2.83822\\
% 10 & 4 & 3.04127 & 3.84138\\
10 & 5 & \textbf{3.55245} & 5.00906\\
10 & 10 & \textbf{6.27067} & 13.818\\ \hline
% 15 & 3 & 3.26562 & 3.87379\\
% 15 & 4 & 3.99006 & 5.34704\\
% 15 & 5 & 4.70297 & 7.09027\\
% 15 & 10 & 8.59396 & 20.0832\\
% 15 & 15 & 12.336 & 39.7061\\ \hline
20 & 3 & \textbf{4.0318} & 4.90157\\
% 20 & 4 & 4.96316 & 6.85085\\
20 & 5 & \textbf{5.89315} & 9.15582\\
20 & 10 & \textbf{10.6386} & 26.3167\\
20 & 15 & \textbf{15.5526} & 52.3773\\
20 & 20 & \textbf{20.4688} & 87.3382\\ \hline
30 & 3 & \textbf{5.49861} & 6.95763\\
% 30 & 4 & 6.83887 & 10.1948\\
30 & 5 & \textbf{8.13738} & 13.3033\\
30 & 10 & \textbf{14.7616} & 38.7993\\
30 & 15 & \textbf{21.4911} & 77.9852\\
30 & 20 & \textbf{28.3982} & 129.726\\
30 & 30 & \textbf{43.0635} & 273.193\\ \hline
50 & 3 & \textbf{8.47438} & 11.0542\\
% 50 & 4 & 10.5455 & 15.8783\\
50 & 5 & \textbf{12.6016} & 21.5876\\
50 & 10 & \textbf{22.9807} & 64.0087\\
50 & 15 & \textbf{33.535} & 128.526\\
50 & 20 & \textbf{44.5519} & 215.052\\
50 & 30 & \textbf{66.6827} & 452.664\\
50 & 50 & \textbf{112.572} & 1196.4 \\ \hline
100 & 3 & \textbf{15.5886} & 21.9397 \\
% 100 & 4 & \textbf{19.3831} & 31.7575 \\
100 & 5 & \textbf{23.1784} & 43.2033 \\
100 & 10 & \textbf{42.3381} & 127.866 \\
100 & 20 & \textbf{80.9126} & 430.52 \\
100 & 30 & \textbf{121.682} & 906.429 \\
100 & 50 & \textbf{204.535} & 2375.7 \\
100 & 100 & \textbf{422.881} & 9126.77 \\
\end{tabular}
\renewcommand{\arraystretch}{1}
\vspace{2ex}
\caption{Running times comparison (in seconds) for Example~\ref{E:Example3}
(higher order derivatives, $d=2$).}
\label{T:Table5}
\vspace{-3ex}
\end{center}
}
\end{table*}

Except for the tests for $n\in\{3,4\}$ and $r=3$, the method using
approach described in Section~\ref{SS:RatDiff-III} was faster than using
the recurrence scheme from Section~\ref{SS:RatDiff-II-Higher}.

\begin{example}\label{E:Example4}
We repeat almost all numerical experiments conducted in 
Examples~\ref{E:Example2} and~\ref{E:Example3} for $d=3$. Results are given in
Table~\ref{T:Table6} and Table~\ref{T:Table7}, respectively.
\end{example}

\begin{table*}[ht!]
{\small
\begin{center}
\renewcommand{\arraystretch}{1.45}
\begin{tabular}{lccccc}
$n$ & $r$ & Floater & accelerated Floater & method from~\S\ref{SS:RatDiff-III} 
& method from~\S\ref{SS:F+S-Der}\\ \hline
2 & 1 & 0.82683 & \textbf{0.76977} & 0.936238 & 0.859879 \\
2 & 2 & 0.855235 & \textbf{0.806501} & 1.08714 & 0.996105 \\ \hline
3 & 1 & 0.876951 & \textbf{0.803918} & 0.989763 & 0.917853 \\
3 & 2 & 0.878117 & \textbf{0.836707} & 1.13669 & 1.08219 \\ \hline
4 & 1 & 0.93218 & \textbf{0.826905} & 1.0242 & 0.957915 \\
4 & 2 & 0.956399 & \textbf{0.907897} & 1.24014 & 1.21526 \\ \hline
5 & 1 & 1.04276 & \textbf{0.876359} & 1.09043 & 1.02799 \\
5 & 2 & 1.06956 & \textbf{0.982089} & 1.33804 & 1.35246 \\ \hline
10 & 1 & 1.835 & \textbf{1.11705} & 1.43842 & 1.37198 \\
10 & 2 & 1.86581 & \textbf{1.32367} & 1.8785 & 1.97703 \\ \hline
20 & 1 & 4.71344 & \textbf{1.5715} & 2.12063 & 2.00936 \\
20 & 2 & 4.75143 & \textbf{2.02847} & 2.95252 & 3.19551 \\ \hline
30 & 1 & 9.36345 & \textbf{2.0241} & 2.7688 & 2.64447 \\
30 & 2 & 9.39308 & \textbf{2.70635} & 3.95168 & 4.41128 \\ \hline
50 & 1 & 23.9363 & \textbf{2.93157} & 4.03636 & 3.90962 \\
50 & 2 & 23.9646 & \textbf{4.06433} & 5.88911 & 6.81092 \\ \hline
100 & 1 & 92.0546 & \textbf{5.22606} & 7.26309 & 7.10418 \\
100 & 2 & 92.1474 & \textbf{7.49552} & 10.8035 & 12.7775 \\ \hline
200 & 1 & 368.895 & \textbf{10.0531} & 14.0438 & 13.869 \\
200 & 2 & 368.967 & \textbf{14.726} & 20.8998 & 25.3485 \\ \hline
300 & 1 & 827.888 & \textbf{15.0003} & 20.9164 & 21.0872 \\
300 & 2 & 826.852 & \textbf{22.0828} & 31.2927 & 39.1476 \\
\end{tabular}
\renewcommand{\arraystretch}{1}
\vspace{2ex}
\caption{Running times comparison (in seconds) for Example~\ref{E:Example4}
(the first two derivatives, $d=3$; cf.~Example~\ref{E:Example2}).}
\label{T:Table6}
\vspace{-3ex}
\end{center}
}
\end{table*}

\begin{table*}[ht!]
{\small
\begin{center}
\renewcommand{\arraystretch}{1.45}
\begin{tabular}{lccc}
$n$ & $r$ & method from~\S\ref{SS:RatDiff-III} 
                                  & method from~\S\ref{SS:RatDiff-II}\\ \hline
3 & 3 & \textbf{1.38741} & 1.40635 \\ \hline
4 & 3 & \textbf{1.54348} & 1.62948 \\
4 & 4 & \textbf{1.71111} & 1.95424  \\ \hline
5 & 3 & \textbf{1.59121} & 1.74595 \\
5 & 4 & \textbf{1.86203} & 2.22166 \\
5 & 5 & \textbf{2.13483} & 2.76566 \\ \hline
10 & 3 & \textbf{2.32797} & 2.74212 \\
% 10 & 4 & \textbf{2.79539} & 3.64522 \\
10 & 5 & \textbf{3.27416} & 4.69733 \\
10 & 10 & \textbf{5.64113} & 12.5317 \\ \hline
20 & 3 & \textbf{3.75143} & 4.72604 \\
% 20 & 4 & \textbf{4.55744} & 6.46642 \\
20 & 5 & \textbf{5.40027} & 8.57394 \\
20 & 10 & \textbf{9.6048} & 24.9434 \\
20 & 20 & \textbf{18.5194} & 84.6039 \\ \hline
30 & 3 & \textbf{5.11798} & 6.63659 \\
% 30 & 4 & \textbf{6.33931} & 9.27392 \\
30 & 5 & \textbf{7.52544} & 12.3671 \\
30 & 10 & \textbf{13.5876} & 35.0617 \\
30 & 20 & \textbf{25.752} & 116.541 \\
30 & 30 & \textbf{38.7294} & 243.808 \\ \hline
50 & 3 & \textbf{7.73327} & 10.4682 \\
% 50 & 4 & \textbf{9.6016} & 14.8763 \\
50 & 5 & \textbf{11.4473} & 20.0135 \\
50 & 10 & \textbf{20.7944} & 57.5285 \\
50 & 20 & \textbf{40.1088} & 192.779 \\
50 & 30 & \textbf{60.8049} & 408.96 \\
50 & 50 & \textbf{101.869} & 1067.28 \\ \hline
100 & 3 & \textbf{14.5509} & 20.4969 \\
% 100 & 4 & \textbf{18.0547} & 29.063 \\
100 & 5 & \textbf{21.5907} & 39.5317 \\
100 & 10 & \textbf{39.4184} & 115.909 \\
100 & 20 & \textbf{76.1972} & 393.102 \\
100 & 30 & \textbf{119.628} & 875.128 \\
100 & 50 & \textbf{202.905} & 2319.34 \\
100 & 100 & \textbf{413.053} & 8899.99 \\
\end{tabular}
\renewcommand{\arraystretch}{1}
\vspace{2ex}
\caption{Running times comparison (in seconds) for Example~\ref{E:Example4}
(higher order derivatives, $d=3$; cf.~Example~\ref{E:Example3}).}
\label{T:Table7}
\vspace{-3ex}
\end{center}
}
\end{table*}

%%%%%%%%%%%%%%%%%%%%%%%%%%%%%%%%%%%%%%%%%%%%%%%%%%%%%%%%%%%%%%%%%%%%%%%%%%%%%%
%%%%%%%%%%%%%%%%%%%%%%%%%%%%%%%%%%%%%%%%%%%%%%%%%%%%%%%%%%%%%%%%%%%%%%%%%%%%%%
\section{Conclusions}                                    \label{S:Conclusions}
%%%%%%%%%%%%%%%%%%%%%%%%%%%%%%%%%%%%%%%%%%%%%%%%%%%%%%%%%%%%%%%%%%%%%%%%%%%%%%
%%%%%%%%%%%%%%%%%%%%%%%%%%%%%%%%%%%%%%%%%%%%%%%%%%%%%%%%%%%%%%%%%%%%%%%%%%%%%%

In this paper, we have presented new algorithms for faster evaluation of 
derivatives of any order of both polynomial and rational B\'{e}zier curves. The 
methods perform well numerically and offer a noticeable boost to the speed of 
the evaluations.

Taking into account the advantages and the disadvantages of the new methods, 
their computational complexity, some technical aspects related to the 
implementation, as well as presented results of the numerical tests, we believe 
that the reader can choose the approach which is best suited for their
particular use case.

Further research is required to rigorously analyse the numerical computation 
errors in the proposed methods. Certainly, one has to first perform the 
numerical analysis of \cite[Alg.~1.1, 2.1 and 2.2]{WCh2020} which are the main 
tools used in our approaches. 

Another interesting issue is to find some bounds for the derivatives of rational 
B\'{e}zier curves (cf., e.g., \cite{Hermann99}, \cite{Shi23} and papers cited 
therein) which follow from the proposed methods. However, it is not so easy 
to give new and useful bounds, because our methods are iterative. This 
problem also needs further research.

%%%%%%%%%%%%%%%%%%%%%%%%%%%%%%%%%%%%%%%%%%%%%%%%%%%%%%%%%%%%%%%%%%%%%%%%%%%%%%
%%%%%%%%%%%%%%%%%%%%%%%%%%%%%%%%%%%%%%%%%%%%%%%%%%%%%%%%%%%%%%%%%%%%%%%%%%%%%%
%\newpage
\bibliographystyle{elsart-num-sort}
\biboptions{sort&compress}
\bibliography{NewBezierDiffEval-rev}

\begin{thebibliography}{10}
\expandafter\ifx\csname url\endcsname\relax
  \def\url#1{\texttt{#1}}\fi
\expandafter\ifx\csname urlprefix\endcsname\relax\def\urlprefix{URL }\fi

\bibitem{BM99}
W.~Boehm, A.~M\"{u}ller, On de {C}asteljau's algorithm, Computer Aided
  Geometric Design 16 (1999) 587--605.

\bibitem{FChPhD}
F.~Chudy, New algorithms for {B}ernstein polynomials, their dual bases, and
  {B}-spline functions, Ph.D. thesis, University of Wroc{\l}aw, available on
  request (2022).

\bibitem{ChW2023}
F.~Chudy, P.~Wo\'{z}ny, Linear-time algorithm for computing the
  {B}ernstein-{B}\'{e}zier coefficients of {B}-spline basis functions, Computer
  Aided-Design 154 (2023) 103434.

\bibitem{DC59}
P.~de~Casteljau, Outillage m\'ethodes calcul (\textit{in French}), Tech. rep.,
  Andr\'e {C}itro\"en {A}utomobile {SA}, Paris (1959).

\bibitem{DC63}
P.~de~Casteljau, Courbes et surfaces \`a p\^{o}les (\textit{in French}), Tech.
  rep., Andr\'e {C}itro\"en {A}utomobile {SA}, Paris (1963).

\bibitem{DC99}
P.~de~Casteljau, De {C}asteljau's autobiography: {M}y time at {C}itro\"en,
  Computer Aided Geometric Design 16 (1999) 583--586.

\bibitem{Farin2002}
G.~Farin, Curves and surfaces for {C}omputer-{A}ided {G}eometric {D}esign. A
  practical guide, 5th ed., Academic Press, Boston, 2002.

\bibitem{Farouki2012}
R.~Farouki, The {B}ernstein polynomial basis: {A} centennial retrospective,
  Computer Aided Geometric Design 29 (2012) 379--419.

\bibitem{Floater91}
M.~Floater, Derivatives of rational {B}\'{e}zier curves, Computer Aided
  Geometric Design 9~(3) (1992) 161--174.

\bibitem{Hermann99}
T.~Hermann, On the derivatives of second and third degree rational {B}\'{e}zier
  curves, Computer Aided Geometric Design 16~(3) (1999) 157--163.

\bibitem{RH21}
A.~Ramanantoanina, K.~Hormann, New shape control tools for rational
  {B}\'{e}zier curve design, Computer Aided Geometric Design 88 (2021) 102003.

\bibitem{RV21}
L.~Romani, A.~Viscardi, Construction and evaluation of pythagorean hodograph
  curves in exponential-polynomial spaces, SIAM Journal on Scientific Computing
  44 (2022) A3515--A3535.

\bibitem{Shi23}
M.~Shi, On the derivatives of rational {B}\'{e}zier curves,
  \url{https://arxiv.org/abs/2303.16156} (2023).

\bibitem{WCh2020}
P.~Wo\'{z}ny, F.~Chudy, Linear-time geometric algorithm for evaluating
  {B}\'{e}zier curves, Computer Aided-Design 118 (2020) 102760.

\end{thebibliography}
%%%%%%%%%%%%%%%%%%%%%%%%%%%%%%%%%%%%%%%%%%%%%%%%%%%%%%%%%%%%%%%%%%%%%%%%%%%%%%
%%%%%%%%%%%%%%%%%%%%%%%%%%%%%%%%%%%%%%%%%%%%%%%%%%%%%%%%%%%%%%%%%%%%%%%%%%%%%%

\end{document}